\input amstex
\documentstyle{amsppt}
\magnification=1200
\NoBlackBoxes
\TagsOnRight
\NoRunningHeads
\vsize=22 truecm
\hsize=16 truecm

%%%%%%%%%%%%%%%%%%%%%%%%%%%%%%%%%%%%%%%%%%%%%%%%%%%%%%%%%%%%%%%%%%%%
\def \real{{\Bbb R}}
\def \complex{{\Bbb C}}

\def \integer{{\Bbb Z}}

\def \trcap{{\cap \negthickspace \negthinspace @!  @! @!| \medspace }}
\def \Id {{\text{Id}\, }}
\def \Fix {{\text{Fix}\, }}
\def \Det {{\text{Det}\, }}
\def \Per {{\text{Per} \, }}
\def \tr {{\text{tr}}}
\def \Int{{\text {int} \, }}

\def\sp{\operatorname{sp}}

\def\AA{{\Cal A}}
\def\BB{{\Cal B}}
\def\CC{{\Cal C}}
\def\DD{{\Cal D}}
\def\EE{{\Cal E}}
\def\FF{{\Cal F}}
\def\GG{{\Cal G}}
\def\HH{{\Cal H}}
\def\II{{\Cal I}}

\def\LL{{\Cal L}}
\def\MM{{\Cal M}}
\def\NN{{\Cal N}}

\def\QQ{{\Cal Q}}
\def\RR{{\Cal R}}
\def\SS{{\Cal S}}
\def\TT{{\Cal T}}
\def\UU{{\Cal U}}

\def\WW{{\Cal W}}

\def\today {\ifcase\month\or January \or February \or March \or
April \or May \or June
\or July \or August \or September \or October \or November \or December
\fi
\number\day~\number\year}

%%%%%%%%%%%%%%%%%%%%%%%%%%%%%%%%%%%%%%%%%%%%%%%%%%%%%%%%%%%%%%%%%%%%

\topmatter

\title\nofrills
Dynamical zeta functions for
analytic surface diffeomorphisms with dominated splitting
\endtitle

\author
Viviane Baladi, Enrique R. Pujals, and Mart\'{\i}n Sambarino
\endauthor

\date
August 2003
\enddate

\address
V. Baladi: CNRS UMR 7586, 
Institut Math\'ematique de Jussieu, 75251 Paris, FRANCE
\endaddress
\email
baladi\@math.jussieu.fr
\endemail

\address
E.R. Pujals:
Instituto de Matem\'atica,
Universidade Federal do Rio de Janeiro, 
CEP 21.945-970, Rio de Janeiro, BRAZIL 
\endaddress
\email
enrique\@impa.br
\endemail

\address
M. Sambarino:
IMERL, Fac. Ingenier\'ia,
Montevideo C.C. 30, URUGUAY
\endaddress
\email
samba\@fing.edu.uy
\endemail

\abstract
We consider a real-analytic compact surface diffeomorphism $f$, for which the
tangent space over the nonwandering set $\Omega$ admits a dominated
splitting. We study  the dynamical determinant
$d_f(z)=\exp-\sum_{n \ge 1} {z^n \over n} \sum_{x \in \Fix^* f^n}
|\Det{ (Df^n(x) - \Id)}|^{-1}$, 
where $\Fix^* f^n$ denotes the set of fixed points
of $f^n$ with no zero Lyapunov exponents.
By combining previous work of Pujals and Sambarino [PS1, PS2]
on $C^2$ surface diffeomorphisms with, on the one
hand,  results of Rugh [Ru1] on hyperbolic analytic
maps, and on the other, our two-dimensional version of the same author's [Ru3]
analysis of one-dimensional analytic dynamics with neutral fixed
points, we prove that $d_f(z)$ is either
an entire function or a holomorphic function in
a (possibly multiply) slit  plane.
\endabstract


\thanks
V.B.  is grateful to PUC-Rio and IMPA for their hospitality at
several stages  of this work and  thanks the
organisers of 2001 l'Odyss\'ee Dynamique in CIRM/Luminy 
and the PRODYN programme of the European Science Foundation
for support.
Visits to Orsay of E.R.P. in 2000 and of E.R.P. and M.S. in 2001
were instrumental in the preparation of this article, we are grateful
to Universit\'e de Paris-Sud and the Brazilian-French 
cooperation program. E.R.P. is partially
supported by FAPERJ-Brazil and the Guggenheim Foundation, and M.S. is
partially supported by the University of Maryland and
CNPq/PRONEX. \endthanks

\endtopmatter

%%%%%%%%%%%%%%%%%%%%%%%%%%%%%%%%%%

\document
\head 1. Introduction 
\endhead

Let $f$ be a real-analytic diffeomorphism 
of a compact  two-dimensional analytic riemannian manifold $M$. 
Our dynamical assumption is that the
tangent space over the nonwandering set $\Omega$ of $f$
admits a dominated splitting, i.e., $T_\Omega M=E \oplus F$, and there are
$C > 0$ and $0 < \lambda < 1$ so that
$$
\| Df^n|_{E(x)}\| \cdot \| Df^{-n}|_{F(f^n(x))}\| \le
C \lambda ^n \, , \forall x \in \Omega\, ,
n \ge 0 \, . \tag{1.1}
$$

We study  the dynamical determinant
$$
d_f(z)=\exp- \sum_{n = 1}^\infty  {z^n \over n} \sum_{x \in \Fix^* f^n}
{1 \over |\Det{( Df^n(x) - \Id)}|} \, ,\tag{1.2}
$$ 
where $\Fix^* f^n$ denotes the (finite) set of hyperbolic fixed points
of $f^n$, i.e., those with no zero Lyapunov exponents. 
(Note that the dominated splitting assumption implies that
each fixed point of $f^n$ has at most one zero Lyapunov exponent.)

In order to state our theorem, we  define 
the {\it $\Omega\setminus P$-isolated periodic points:} they are the elements
$p$ of the set $\Per$ of  periodic points
of $f$ which are not in  $\overline{ (\Omega \setminus \Per)}$,
i.e., which admit a neighbourhood $U_p$ with $U_p \cap \Omega \subset \Per$.
(Besides isolated hyperbolic sinks, sources, and saddles, this set only
contains periodic points which do not contribute to the chaotic dynamics.)
In Section 4, we shall recall the  decomposition of $\Omega$ into
periodic, quasi-periodic, and almost hyperbolic components
from [PS2]. This will be our starting point in the analysis of $f$.
We just mention here the fact that the set  $\NN$ of nonhyperbolic
periodic points of $f$ which are {\it not} $\Omega\setminus P$-isolated periodic points is 
(empty or) finite.
If $p\in \Per$, we write  $P=P(p)\ge 1$ for its minimal
period and $\lambda_E=\lambda_E(p)$, $\lambda_F=\lambda_F(p)$ for its multipliers, i.e., eigenvalues of
$Df^P(p)$ with $|\lambda_E| < |\lambda_F|$  (both multipliers are real because of the
dominated splitting). 
We associate to each $p \in \NN$  the following subset of $\complex$:
$$
\Sigma(p)=
\cases
\{ z\mid  z^P \in  [-1, 1]  \} \, ,
&\hbox{if $\lambda_F=-1$ and $|\lambda_E| < 1$,}\cr
\{  z \mid z^P \in  [\min (0, \lambda_E) , 1]  \} \, ,
&\hbox{if $\lambda_F=+1$ and $|\lambda_E| < 1$,}\cr
\{  z \mid z^P \in  |\lambda_F|^{-1}\cdot  [- 1, 1] \}  \, ,
&\hbox{if $\lambda_E=-1$ and $|\lambda_F|>1$,} \cr
\{  z \mid z^P \in |\lambda_F|^{-1}\cdot [ \min(0, \lambda_F^{-1}) ,1 ] \} \, ,
&\hbox{if $\lambda_E=+1$ and $|\lambda_F|>1$.} \cr
\endcases \tag{1.3}
$$

Our main result can now be summarized as follows (see Section~ 4 for
more):

\proclaim{Theorem A}
Let $f : M \to M$ be a real-analytic diffeomorphism of
a compact analytic riemannian surface. Assume that
$f$ admits a dominated splitting \thetag{1.1}  over its nonwandering
set $\Omega$. 
Let $\NN$ be the empty or finite set of nonhyperbolic periodic
points which are {\it not} $\Omega\setminus P$-isolated periodic points.
Then $d_f(z)$ is is holomorphic and nonzero
in the open unit disc and admits a holomorphic extension to the plane,
slit plane, or multiply slit plane defined by
$$
\{ z \in \complex \mid  1/z \notin \cup_{p\in \NN}\Sigma(p) \}\, .\tag{1.4}
$$
\endproclaim

We conjecture that the endpoints of the slits are nonpolar singularities.
It is an open question whether $d(z)$ may be analytically
continued across the open slits (to different sheets of a Riemann surface).
See below for more involved conjectures and questions.

\smallskip
Theorem A immediately implies:
\roster
\item
If $\NN$ is empty, i.e., if all nonhyperbolic points
are $\Omega\setminus P$ isolated, then $d_f(z)$  is an entire function with no zeroes 
in the open unit disc.  
\item
If there exist points in $\NN$ with 
$|\lambda_F|=1$, then $d_f(z)$ is analytic and
nonzero in the disc of radius $1$, with a possibly
nonpolar singularity at $z= 1$ or $-1$, and it
admits an analytic extension to a (possibly multiply) slit plane. 
\item
If there exist  points in $\NN$ with $|\lambda_E| =1$, but no points
in $\NN$ with $|\lambda_F|=1$,
letting $|\lambda_F|^{1/P}$ be the smallest modulus of
$P$-th roots of $F$-multipliers in $\NN$,
 then $d_f(z)$ is analytic and nonzero in the unit disc, and it may be analytically
extended to the disc of radius $|\lambda_F|^{1/P}>1$, with (finitely many) possibly non polar
singularities on its boundary, and a further analytic extension
to a (possibly multiply)  slit plane.
\endroster

\smallskip

We
next say a few words about the proof of Theorem ~A, sketching the contents of the paper.
If $\NN$ is empty, we shall see in Section~4 that $f$ is uniformly hyperbolic
on a compact invariant subset $\Lambda$  of its wandering set which
contains all the nonisolated hyperbolic periodic points.
The results of Rugh [Ru1, Ru2] on the dynamical determinants
of hyperbolic analytic maps immediately
imply that $d_f(z)$ is an entire function.
The key point in Rugh's analysis, inspired by
Ruelle's [Rue1] seminal study (Ruelle only considered the case when the dynamical
foliations are analytic),  was to express
$d_f(z)$ as a quotient of the Grothendieck-Fredholm
determinants of two nuclear operators, proving also that
zeros in the denominator are always cancelled by the
numerator. The nonanalyticity of the dynamical foliations
can be disregarded by working with two contracting and
analytic half-inverses of $f$, in appropriate coordinates.
In practice, Rugh [Ru2] constructs a symbolic model
for a real-analytic hyperbolic map, starting from a Markov partition.

If $\NN$ is not empty we must modify Rugh's model
to investigate $d_f(z)$. The description of the
corresponding {\it almost hyperbolic  real-analytic}
symbolic model $\hat f$
is carried out in Sections ~2 and 3, while Section ~4
discusses how to reduce from our surface diffeomorphism $f$ to  $\hat f$.
In a nutshell, we discuss in Section ~4 
Markov partitions for $f$, describing how they contain both
``good'' (i.e., of  hyperbolic type) and
``bad'' rectangles (those which contain an element of $\NN$).
The dynamical determinant
$d_f(z)$ is morally the (regularised) determinant of a transfer operator
$\widehat \LL$
analysed in Sections 2--3.
The building blocks of $\widehat \LL$  are either ``good,'' 
and of the type studied in [Ru1, Ru2], or ``bad'' and approximate direct products
of a one-dimensional hyperbolic operator with a one-dimensional
parabolic operator,  studied in another work of Rugh [Ru3]. 
For the parabolic operator,
we use a normal form [H] to adapt the analysis of
one-dimensional analytic dynamics with neutral fixed
points in  [Ru3], to our setting. 
(In [Ru3], the nondiscrete spectrum of the operator was
the compact interval $[0,1]$.)

More precisely, we describe in Section~ 2 the almost hyperbolic model
and introduce the building blocks $\widehat \LL_{kj}$ of the symbolic transfer operators
as well as the Banach spaces $\BB_k$ they act on.
In Section~3.A, we analyse the spectrum of the ``bad''
$\widehat \LL_{kk}$.
The crucial tool to do
this is an approximate  Fatou coordinate for parabolic points.
Section~3.B contains a complete description of the spectrum
and of the regularised determinant of the symbolic model.
We combine  Sections 2--3 wih Section 4 in 
Section 5.A, using a sequence of Markov partitions with
diameter going to zero, to show that, for every
neighbourhood of the ``slit plane'' \thetag{1.4} in Theorem~A, $d_f(z)$ is  analytic
outside of this neighbourhood.

\smallskip

To keep the paper reasonably short, we do not reproduce
the arguments of Rugh when they can be used without nontrivial
modifications. Note also that for the sake of simplicity,
in Sections 2--4, we (mostly) restrict to the case where all elements
of $\NN$ are {\it parabolic fixed} points in a strict sense (i.e., the order of
$f-\Id$ is equal to two) with a nonhyperbolic multiplier equal
to $+1$ (and not $-1$). The reduction from the general case to
this setting is explained in Section ~5.B.

The appendix is devoted to the construction of adapted metrics
in our setting.

\medskip

Finally, in order to state a conjecture motivated by our result,  we
recall some definitions.
A  $u$-Gibbs state is an ergodic invariant
probability measure, whose induced measures along the
Pesin unstable manifolds are absolutely continuous with respect
to Lebesgue (see [PeS]).  An invariant
probability  measure $\mu$ is called a physical measure if there
is a set of positive Lebesgue measure of points $x$ 
so that ${1\over n} \sum_{k=0}^{n-1} \delta_{f^k(x)}$
weak-$*$ converges to $\mu$ as $n \to \infty$.
We say that a compact invariant set $\Lambda$ is an attracting set if
there is an open neigbourhood $\widehat \Lambda$ such that
$\Lambda =\cap_{n\ge 0} f^n(\widehat \Lambda)$.
Let $\Lambda_j$ be a basic set  in the 
decomposition of $\Omega$ from [PS1, PS2] detailed in Section~ 4,
and assume that $\Lambda_j$ is
an attracting set which does not contain any periodic
point with $|\lambda_F|=1$.
Using the information in the appendix below~---~
in particular Lemma A.5 for $F$-hyperbolic points
which gives the ``mostly contracting'' condition;
density of the strong unstable leaves follows from the fact that each
mixing basic subset is a homoclinic class
$\overline{W^s(p) \trcap W^u(p)}$~---~
the results of Bonatti-Viana [BoVi] may be adapted to
$f|\Lambda_j$, proving that
it enjoys a single $u$-Gibbs state  which is also a physical measure.
Let us call SRB measure a $u$-Gibbs
state which is also a physical measure. (In particular, the
Dirac mass at a hyperbolic sink is an SRB measure.)

By ``exponential rate of mixing'' for an $f$-invariant probability
measure $\mu$, we mean that there is $\tau < 1$ so that
the correlation function satisfies
$$
|\rho_{\varphi, \psi}(k)|=|\int \varphi \circ f^k \psi \, d\mu -\int \varphi \, d\mu
\int \psi \, d\mu | \le C_{\varphi \psi} \tau^{|k|}\, , \tag{1.5}
$$
for all $k \in \integer$ and all
Lipschitz $\varphi$, $\psi$, with $C_{\varphi \psi}$ 
depending on the Lipschitz norms of $\varphi$ and $\psi$.
Let $\Lambda_{i_j}$ be a (topologically mixing for an iterate $f^{n_i}$)  basic subset
(from the  decomposition of $\Omega$  recalled in Section~4),  which 
is attracting and does not contain any periodic
point with $|\lambda_F|=1$.
The results of Castro [Cas] (see also Dolgopyat [Do]) indicate that
the unique SRB measure $\mu$ for $f^{n_i}$ on $\Lambda_{i_j}$ furnished by Bonatti-Viana [BoVi] 
has exponential rates of mixing
(for Lipschitz observables). 
Define the {\it
analytic correlation spectrum} of  such an attracting basic subset $\Lambda_{i_j}$
and its SRB measure $\mu$ to be the union over all
pairs of analytic observables $\varphi$, $\psi$ (extending holomorphically on a fixed complex neighbourhood
of $M$) of the singular set of the Fourier transform of the correlation  function:
$$
\hat \rho_{\varphi \psi}(\omega)=
\sum_{k \in \Bbb Z} e^{i \omega k} \rho_{\varphi \psi} (n_i k) \, .
$$  
Exponential decay of the correlation function implies that
$\hat \rho$ is analytic in the strip $|\Im \omega|< \log 1/\tau$.

If $\NN=\emptyset$, the order $D\ge 0$ of $1$ as a zero of $d_f(z)$ coincides
with the number of  SRB measures of $f$. 
(As observed above, we are in a hyperbolic situation.
The statement for the multiplicity of the pole at $1$ of the dynamical
zeta function $\zeta_f(z)=\exp \sum_{n\ge 1}{z^n\over n }
\sum_{x\in \Fix f^n} \lambda_F(x)^{-n/P(x)}$,
weighted with $\lambda_F^{-1/P}$, follows e.g. from
[Rue2]. Then use that $d_f(z) \zeta_f(z)$ is holomorphic and nonvanishing
in a disc of radius larger than $1$.)

\proclaim{Conjecture B}
Assume that $\min_{\NN} |\lambda_F| > 1$.
We conjecture that the the order $D\ge 0$ of $1$ as a zero of $d_f(z)$ coincides
with the number of (ergodic) SRB measures of $f$ (which coincides
with the number of attracting basic sets added to the number of hyperbolic
sinks). 
Furthermore, under the same assumption,  let $\mu$ be
the SRB measure of an attracting basic set $\Lambda_j$ in $\Lambda$.
We  conjecture that $\omega\ne 0$  is in the analytic correlation spectrum  only if:
either $d_{f|{\Lambda_j}}(z)$ is not holomorphic at $z=\exp (i \omega)$
or it vanishes there. We ask whether the other implication holds.
In particular,  we conjecture that the SRB measure $\mu$
has exponential rates of mixing  (if and only if $\Lambda_j$ is mixing),
if and only if  $z=1$ is the only zero
of modulus one of $d_{f|_{\Lambda_j}}(z)$.
\endproclaim

Conjecture B states in particular that
the presence of a ``gap'' in the dynamical determinant  
(of a transitive component) reflects 
exponential mixing of the SRB measure in the setting of (analytic) surface diffeomorphisms
enjoying dominated splitting. The only setting where we
know a (proved) analogue of this statement is $S$-unimodal interval
maps [Ke].

{\bf Acknowledgements:}
V.B. is grateful to  Artur Avila, S\'ebastien Gou\"ezel, and
Hans Henrik Rugh for useful comments. 
\bigskip

\head 2. The symbolic maps and their transfer operators
\endhead

\subhead 2.A  Almost hyperbolic analytic maps
\endsubhead
The key is to reduce (using suitable coordinate charts
on Markov covers close
to small enough Markov partitions) the problem to a
variant of the symbolic model introduced in [Ru1]: (real)-analytic
hyperbolic maps.   We shall call ``almost hyperbolic analytic maps''
our variant, where some of the
building blocks are associated to periodic points with a
neutral multiplier.
It is convenient to use the following
open ``petals'' in $\complex$, associated
to  real numbers $r_0>0$, $\theta_0\in (0, \pi)$ by
$$
\UU(\theta_0,r_0)=
\{r e^{i\theta} \in \complex \mid 0<r < r_0\, ,
\,  -\theta_0 < \theta < \theta_0 \}\, .\tag{2.1}
$$

\remark{Remark 2.1 (Fixed points or periodic points)}
The symbolic model of this section is adapted to the situation
where all nonhyperbolic periodic points are {\it fixed} points with
an eigenvalue $+1$ (but no eigenvalue $-1$).
It is not difficult, although cumbersome, to lift this restriction, we shall
do this  Section 5.B.
\endremark

\definition {Definition 2.2
(The model: almost hyperbolic analytic surface map)}
An almost hyperbolic analytic surface map $\hat f$ consists in a finite set
$\SS = \SS_0 \cup \SS_1$, with $\SS_0$, $\SS_1$ nonempty and
disjoint, and data
$$
\eqalign{
&
[(\DD^1_i\times \DD^2_i\subset \complex \times \complex )\, ,
(t_{ij}\in \{0,1\})\, ,
(\hat f_{ij}: \DD^1_i \times \DD^2_i \to \DD^1_j\times \complex\, , \forall \, t_{ij}=1)\mid i, j \in \SS]
\, , 
}
$$
satisfying the following assumptions:
The $\SS \times \SS$ matrix $(t_{ij})_{i,j\in I}$ is irreducible
with no wandering states, and $t_{ij}=1$ with $(i,j)\in \SS_0\times \SS_0$
if and only if $i=j$.
Each $\hat f_{ij}$  admits a real-analytic 
extension $(\hat f_{ij}^{1}, \hat f_{ij}^{2})$ to a neighbourhood of 
$\DD^1_i \times \DD^2_i$ and 

(H) For all $i\in  \SS_1$, $\DD^1_i$ and $\DD^2_i$ are connected compact
subsets of the complex plane with  $C^1$  boundaries.
If $i$ or $j \in \SS_1$ then:
\smallskip

\noindent (H.1) $\hat f_{ij}$ is contracting in the first coordinate, i.e.,
$
\hat f_{ij}^{1}({\DD^1_i}\times{\DD^2_i})
\subset \hbox{Int} (\DD^1_j)$.

\noindent (H.2) $\hat f_{ij}$ is expanding in the second coordinate,
i.e., there is a real-analytic function 
$\phi_{ij,s}$ (the ``partial inverse'')
defined on a neighbourhood of 
${\DD^1_i}\times{\DD^2_j}$,
such  that
for each $(w_1, z_2) \in {\DD^1_i}\times {\DD^2_j}$
the image $w_2=\phi_{ij,s} (w_1, z_2)$ lies in the interior
of $\DD^2_i$ and is the unique solution  of
$$
\hat f^{2}_{ij}(w_1, w_2) = z_2 \, .\tag{2.2}
$$

\smallskip

(P) If $i=j\in \SS_0$   there are an integer $\nu_i \ge 1$,  
and real numbers $\pi/(2\nu_i)<\theta_i< \tilde \theta_{i} < \pi/\nu_i $
and $\tilde r_i > r_i>0$ so that $\hat f_{ii}$
is described in one of the two   ``partially hyperbolic''
forms:

\smallskip
\noindent (P.a) $\DD^1_i$ is a closed disc centered
at the origin. The boundary of the compact connected simply connected set 
$\DD^2_i$ is $C^1$ except at one point
which is assumed to be the origin.
$\DD^2_i$ contains
$\UU(\theta_i, r_i)$ and is contained in the closure of
$\UU(\tilde \theta_i, \tilde r_i)$. The map $\hat f_{ii}$
fixes the origin and is contracting in the first coordinate
in the sense of (H.1). There is a real-analytic 
function  $\phi_{ii,s}(w_1,z_2)$  defined on a neighbourhood of 
$\DD^1_i \times  \DD^2_i$,
which is the unique solution  of
$\hat f^{2}_{ii}(w_1, \cdot) = z_2$ there, and which has
the normal form: 
$$
\phi_{ii,s} (w_1, z_2)= z_2 -  z_2^{1+\nu_i }+z_2^{2+\nu_i} 
\tilde \phi_{ii,s}(w_1,z_2) \, . \tag{2.3}
$$

\noindent (P.b)   $\DD^2_i$ is a closed disc centered
at the origin. The boundary of the compact connected simply connected set
$\DD^1_i$ is $C^1$ except at one point
which is assumed to be the origin. $\DD^1_i$ contains $\UU(\theta_i, r_i)$
and is contained in the closure of
$\UU(\tilde \theta_i, \tilde r_i)$.
The map $\hat f_{ii}$  fixes the origin and is expanding in the second coordinate
in the sense of (H.2).
The (real-analytic) map $\phi_{ii,u}(w_1,z_2)=
\hat f_{ii}^1 (w_1, \phi_{ii,s}(w_1,z_2))$ enjoys the normal form
$$
\phi_{ii,u}(w_1,z_2) = w_1 -   w_1^{1+\nu_i }+  
w_1^{2+\nu_i} \tilde \phi_{ii,u}(w_1,z_2)\, .\tag{2.4}
$$
 \enddefinition

\remark{Remark 2.3} 
\roster
\item {\bf (Admissible sequences -- $\II_j^i$)} A symbol sequence $\vec{\imath}=(i_1, \ldots, i_n) \in \SS^n$ is
called admissible if
$t_{{i_k}, {i_{k+1}}} =1$ for every  $k=1, \ldots, n-1$.
Write 
$
\II^j_i = \real \cap \DD^j_i$ for $j=1, 2$, $i\in  \SS$.
\item {\bf (Attracting and repelling petals)}
In case (P.a), up to shrinking the domains,
the normal form implies that for each $w_1\in \DD_i^1$ the map $\phi_{ii,s}(w_1,\cdot)$ 
sends the ``attracting petal'' $\DD^2_i$ injectively into 
$\Int \DD^2_i \cup \{0\}$. In case (P.b),
for each $z_2\in \DD_i^2$ the map $\phi_{ii,u}(\cdot ,z_2)$  sends the ``repelling
petal''
$\DD^1_i$ injectively into  $\Int \DD^1_i \cup \{0\}$. 

\item {\bf (Invertibility)} 
If $|\det D\hat f_{ij}|$ does not vanish
on ${\II^1_i}\times  \phi_{ij,s}({\II^1_i}\times {\II^2_j})$
we say that $\hat f$ is invertible. 
To fix ideas, and since this will be the case
in our application to surface diffeomorphisms, we assume from now on that $\hat f$
is invertible.
\endroster
\endremark

\medskip
\noindent {\bf Pinning coordinates (half-inverses)}
\smallskip

To define the transfer operator, Rugh introduces  ``pinning coordinates:''

\definition{Definition 2.4 (Pinning coordinates/half inverses)}
Let $\hat f$ be almost hyperbolic analytic.
For each pair $(i,j)$ with $t_{ij}=1$, the {\it pinning coordinates}
are the two real-analytic
maps $\phi_{ij,s}$, $\phi_{ij, u}$ defined  on 
a neighbourhood of ${\DD^1_{i}} \times  {\DD^2_{j}}$ by
$$
\cases
\phi_{ij,s} (w_1, z_2) \in \DD^2_i \, ,\cr
\phi_{ij,u}(w_1, z_2) = \hat f^{1}_{ij} (w_1, \phi_{ij,s}(w_1,z_2)) \in \DD^1_j \, .
\endcases
\tag{2.5}
$$
In the hyperbolic case (H),
the pinning coordinate maps  ${\DD^2_{j}}$  
for each $w_1\in {\DD^1_{i}}$, respectively ${\DD^1_{i}}$  
for each $z_2\in {\DD^2_{j}}$,
  injectively into the interior
of $\DD^2_{i}$, respectively $\DD^1_{j}$.
In the parabolic case (P.a) 
$\phi_{ii,u}(\cdot, z_2) $ maps  ${\DD^1_{i}}$ injectively into the interior of $\DD^1_{i}$,
while $\phi_{ii,s} (w_1, \cdot)$ maps $\DD^2_i$ injectively into  $\Int \DD^2_i \cup \{0\}$.
Case (P.b) has the natural symmetric characteristics.
\enddefinition

It is easy to check that Definition~2.4 implies 
$$
\hat f_{ij}(w_1, \phi_{ij,s}(w_1, z_2))=
(\phi_{ij,u}(w_1, z_2), z_2)\, , \quad \forall (w_1,z_2) \in \DD^1_i \times \DD^2_j\, .\tag{2.6}
$$
Proposition 2.6 will show that pinning coordinates  allow
to ``pin down'' the whole orbit up to time $n$
by knowing only the first coordinate of the
initial position and the second coordinate of
the final position. They are  in some sense
``half-inverses'' for the map. They have Fried's [Fr]
``cross maps'' as an avatar. See also the ``implicit (or
macroscopic) coordinates''
used by Palis-Yoccoz [PY].

A word about terminology: By definition of an almost hyperbolic
analytic map, each $\phi_{ij, s}(w_1,\cdot)$
is a diffeomorphism onto its image for fixed $w_1$, 
so that $\phi_{ij, s}$ deserves to be called a ``coordinate.''
If $\hat f$ is invertible, then  each $\phi_{ij, u}(\cdot,z_2)$
is a diffeomorphism onto its image for fixed $z_2$  and both pinning maps deserve
to be called ``coordinates.'' 
Finally, the notation $\phi_{ij, u}$ for the (weakly) contracting
direction of $\hat f_{ij}$ and $\phi_{ij, s}$ 
for its (weakly) expanding direction is used, not only for the
sake of compatibility with published literature, but
because the transfer operator 
is (at least morally, see Remark~3.6)
associated to the inverse dynamics, which exchanges stable and
unstable directions. This remark also applies to the ``attracting'' vs.
``repelling petal'' terminology (the adjectives refer to
the inverse map) in Remark 2.3 and the following definition:

\definition{Definition 2.5 (Complements to Definition 2.2)}
\roster 
\item {\bf (Convenient extensions)}
The reader is invited to check (Schwarz lemma) that
Definition 2.2 implies (up to slightly changing
the domains $\DD_{i}^{1,2}$) that there are compact connected $\widetilde \DD_i^{1,2}$ with
$\DD_i^1\subset \Int \, \widetilde \DD_i^1$
and $\DD_i^2\subset \Int \, \widetilde \DD_i^2$
so that, if $(i,j)\not \in \SS_0\times \SS_0$
then  
\itemitem {(1.a)} $\hat f_{ij}$ may be extended analytically to
$\widetilde \DD_i^1 \times \widetilde \DD_i^2$,
with $
\hat f_{ij}^{1}(\widetilde {\DD^1_i}\times\widetilde {\DD^2_i})
\subset \Int (\DD^1_j)$ in (H.1);
\itemitem {(1.b)} $\phi_{ij, s}$ is defined on 
$\widetilde {\DD^1_i}\times\widetilde {\DD^2_j}$, and
for each $(w_1, z_2)$ there,
$w_2=\phi_{ij,s} (w_1, z_2)$ lies in the interior
of $\DD^2_i$ 
in (H.2);
\itemitem{(1.c)} if $i=j \in \SS_0$,  the hyperbolic conditions of (P.a), (P.b)
hold for the extended domain, and
we may assume that $\phi_{ii,s}(\DD_i^1\times\DD_i^2)
\subset \Int \widetilde \DD_i^2$ in case (P.a) and
$\phi_{ii,u}(\DD_i^1\times\DD_i^2)
\subset \Int \widetilde \DD_i^1$  in case (P.b).

\item {\bf (Case (P.b) -- Attracting petals)}
In case (P.b),  it follows from the assumptions
(up to slightly changing $\DD_i^{1,2}$, $\widetilde \DD_i^{1,2}$) that
for each $z_2\in \DD^2_i$, the map $\phi_{ii,u}(\cdot,z_2)$ sends an ``attracting
petal''   $\DD^{1,-}_i:=e^{i \pi/\nu_i}\DD^1_{i}$ injectively on a domain
containing $\DD^{1,-}_i$, and the inverse
transformation $\phi_{ii,u}^{-1}(\cdot,z_2)$
maps  $\DD^{1,-}_i$ into $\DD^{1,-}_i \cup \{0\}$.  Up to taking a smaller attracting
petal,
we can assume additionally that $\phi_{ij,u} (\DD^1_i, \DD^2_j) \cap\DD^{1,-}_j=\emptyset$
for all $j\in \SS_1$ with $t_{ji}\ne 0$. 
\endroster
\enddefinition

From now on we  assume that $\nu_i= 1$ if $i \in \SS_0$, i.e.,
that all nonhyperbolic transitions are parabolic in the strict
sense. We write $-\DD^1_i$ for the attracting petal
$\DD^{1,-}_i$. (We shall explain in Section~5.B how to reduce to
this situation from the general case.) 

\smallskip

The following result is the key to iterating the ``half-inverses.''
It hinges heavily on the analytic structure.

\proclaim {Proposition 2.6 (Iterating pinning coordinates in $\SS^{n+1}\setminus \SS_0^{n+1}$)}
For each $n \ge 1$ and any admissible symbol
sequence $\vec{\imath} \in \SS^{n+1}\setminus \SS_0^{n+1}$, there are uniquely
defined   iterated pinning maps, real-analytic in a neighbourhood
of  $\DD^1_{i_1} \times \DD^2_{i_{n+1}}$ 
$$
\phi_{\vec{\imath},s}^{(n)}: \DD^1_{i_1} \times \DD^2_{i_{n+1}} \to 
\DD^2_{i_{1}}\, , \,\quad
\phi_{\vec{\imath},u}^{(n)}: \DD^1_{i_1} \times \DD^2_{i_{n+1}} \to 
\DD^1_{i_{n+1}}\, , \tag{2.7}
$$
mapping   $\DD^1_{i_1} \times \DD^2_{i_{n+1}}$
into the interior of $\DD^2_{i_{1}}$,
respectively $\DD^1_{i_{n+1}}$, in such a way that
$$
\eqalign{
\hat f^{(n)}_{\vec \imath} 
(w_1, \phi^{(n)}_{\vec\imath , s}(w_1, z_2))
&:= \hat f_{i_n i_{n+1}} \circ \cdots 
\circ \hat f_{i_1 i_2} (w_1,\phi^{(n)}_{\vec\imath , s}(w_1, z_2)) \cr
&=(\phi^{(n)}_{\vec \imath , u}(w_1, z_2), z_2) 
\, \quad \forall (w_1, w_2) \in \DD^1_{i_1} \times \DD^2_{i_{n+1}}\, .\cr}\tag{2.8}
$$
\endproclaim

\demo{Proof of Proposition 2.6}
If there are no consecutive symbols in
$\SS_0$ in $\vec \imath$, case (P) never occurs, and we
are in the setting of [Ru1]. Let us recall his
proof for completeness: set $\phi^{(1)}_{i_1 i_2, u}=\phi_{i_1 i_2, u}$
and $\phi^{(1)}_{i_1 i_2,s}=\phi_{i_1 i_2,s}$. For $n \ge 2$,
assume by induction that the maps
$$
\cases
\phi^{(1)}_{i_n i_{n+1}, s} : \DD^1_{i_n}  \times \DD^2_{i_{n+1}} 
\to \DD^2_{i_n} \, , \cr
\phi^{(n-1)}_{i_1 \ldots i_n, u} : \DD^1_{i_1} \times \DD^2_{i_n}
\to \DD^1_{i_n} \, , \cr
\endcases
$$
have been defined and satisfy the required properties.
Then, for each fixed  $(w_1,z_2)$ in
$\DD^1_{i_1}\times \DD^2_{i_{n+1}}$,  one easily
checks that the (real-analytic) map
$$
\Phi_{w_1,z_2}(\xi_1,\xi_2) =
(\phi^{(n-1)}_{i_1 \ldots i_n, u}(w_1, \xi_2), 
\phi^{(1)}_{i_n i_{n+1}, s}(\xi_1,z_2)) \tag{2.9}
$$
is contracting in the sense that
$
\Phi_{w_1,z_2}(
{\DD^1_{i_n}} \times {\DD^2_{i_n}}) \subset
\hbox{Int} (\DD^1_{i_n} \times \DD^2_{i_n} )
$.
One deduces from this (Lemma~1 in [Ru1], see also [Rue1]) that $\Phi_{w_1,z_2}$
possesses a unique fixed point $(\xi_1^*, \xi_2^*)
\in \hbox{Int} (\DD^1_{i_n} \times \DD^2_{i_n} )$,
which depends analytically on $w_1$ and $z_2$. Finally,
since $\phi^{(n-1)}_{i_1 \ldots i_n, s}$ and $\phi^{(1)}_{i_n i_{n+1}, u}$ 
also exist by induction,
define the pinning maps by
$$
\cases
\phi^{(n)}_{\vec\imath, u} (w_1, z_2) &= 
\phi_{i_n i_{n+1}, u}^{(1)}(\xi_1^*(w_1,z_2), z_2)\, , \cr
\phi^{(n)}_{\vec\imath, s} (w_1, z_2) &= 
\phi^{{(n-1)}}_{i_1 \ldots i_n, s} (w_1, \xi_2^*(w_1, z_2)) \, .
\endcases\tag{2.10}
$$
Indeed, the induction assumption together
with the fixed point property imply that
$$
\eqalign
{
\hat f_{i_n i_{n+1}}\hat f^{(n-1)}_{i_1 \ldots i_n} 
(w_1, \phi^{(n-1)}_{i_1 \ldots i_n, s} (w_1, \xi_2^*))
&=\hat f_{i_n i_{n+1}}(\phi^{(n-1)}_{i_1 \ldots i_n, u} (w_1, \xi_2^*), \xi_2^*)\, ,\cr
&=\hat f_{i_n i_{n+1}}(\xi_1^*, \phi^{(1)}_{i_n i_{n+1}, s}(\xi_1^*, z_2))\cr
&= ( \phi^{(1)}_{i_n i_{n+1}, u}(\xi^*_1 ,z_2),z_2) \, . \cr
}
$$
Uniqueness follows by induction and uniqueness of $\xi_1^*$, $\xi_2^*$.

In fact, the above argument shows that the iterated pinning coordinates
map   $\widetilde \DD^1_{i_1} \times \widetilde \DD^2_{i_{n+1}}$
into the interior of $\DD^2_{i_{1}}$,
respectively $\DD^1_{i_{n+1}}$.

Next,  assume for the moment that there is
a single occurrence of $i_T=i_{T+1}$ in 
$\SS_0$, at time $T \in \{1, \ldots, n\}$. 
The structure of the induction means we only need
to consider $T=1$ or $n$: There are thus
four possibilities and we consider first the case
$T=n$ and type (P.a). Our starting point is then the pair
$$
\phi^{(1)}_{i_n i_{n+1}, s} : \DD^1_{i_n}  \times \DD^2_{i_{n+1}} 
\to \widetilde \DD^2_{i_n} \, , \qquad
\phi^{(n-1)}_{i_1 \ldots i_n, u} : \DD^1_{i_1} \times \DD^2_{i_n}
\to \DD^1_{i_n}\, . 
$$
For each $(w_1,z_2)$ in
$\DD^1_{i_1}\times \DD^2_{i_{n+1}}$,  the map \thetag{2.9}
is contracting in the sense that
$$
\Phi_{w_1,z_2}(
{\DD^1_{i_n}} \times {\widetilde \DD^2_{i_n}}) \subset
\hbox{Int} (\DD^1_{i_n} \times \widetilde \DD^2_{i_n} )\, .
$$
Therefore, $\Phi_{w_1,z_2}$
possesses a unique fixed point $(\xi_1^*, \xi_2^*)
\in \DD^1_{i_n}\times \widetilde \DD^2_{i_n}$, 
depending analytically on $w_1$ and $z_2$.
Using the inductive assumption
(see also the remark at the end of the hyperbolic
argument) we may again define the pinning maps by
\thetag{2.10}.
The case where $T=1=n-1$ and (P.b) hold is similar.

If $T=n$ and (P.b), or $T=1=n-1$ and (P.a), hold 
then we need to use the inclusion
from (H.1), respectively (H.2), at $i_{n-1} i_n$.

If there are never more than two consecutive symbols in
$\SS_0$ but possibly several such pairs (separated by symbols
in $\SS_1$), the argument just described also applies.

Finally, the case where there are (possibly more than
one group of) three or more consecutive occurrences of symbols in
$\SS_0$ may be dealt with by an easy induction on the number of
such occurrences (using again the property at the
end of the first part of the proof). The idea is to first
consider an $\SS_0 \SS_0 \SS_1$ (or $\SS_1 \SS_0 \SS_0$)
event (as explained above),
and then add one by one the preceding (or following)
symbols in $\SS_0$.
\qed
\enddemo

\remark{Remark 2.7 (Pinning coordinates as coordinates)}
For each $n \ge 1$ and
each admissible $\vec{\imath} \in \SS^{n+1}\setminus \SS_0^{n+1}$,
Proposition~ 2.6
produces a parametrisation of the subset $\DD_{\vec \imath}$ of 
$\DD^1_{i_1} \times \DD^2_{i_1}$ consisting of those
points $(w_1, w_2)$ such that $\hat f_{i_n i_{n+1}} \circ \cdots
\circ\hat f_{i_1 i_2}(w_1, w_2)$ is well-defined and belongs to $\DD^1_{i_{n+1}} \times 
\DD^2_{i_{n+1}}$.
More precisely, $\DD_{\vec \imath}$ is the isomorphic 
image of $\DD^1_{i_1} \times \DD^2_{i_{n+1}}$ 
under the transformation
$
(w_1, z_2) \mapsto (w_1, \phi^{(n)}_{\vec \imath, s} (w_1, z_2)) 
$.
(Real-analyticity implies that this transformation maps  
$\II^1_{i_1}\times
\II^2_{i_{n+1}}$ into $\II^1_{i_1} \times \II^2_{i_1}
\subset \real \times \real$.)
\endremark

\smallskip

\noindent An important consequence of Proposition~2.6 is the
following statement (see also [Ru1]):

\proclaim{Corollary 2.8 (Iterating  hyperbolic analytic maps)}
For each $n\ge 1$ and every admissible symbol
sequence $\vec{\imath} \in  \SS^{n+1}\setminus \SS_0^{n+1}$,
there exists an almost hyperbolic analytic map 
$\hat f^{(n)}_{\vec {\imath}}
:\DD_{i_1 \ldots i_{n+1}} \to \DD^1_{i_{n+1}}\times \DD^2_{i_{n+1}}$  such that
$$
\hat f^{(n)}_{\vec {\imath}} |_{\DD_{i_1 \ldots i_{n+1}} }=
\hat f_{i_{n-1} i_n} \circ \cdots \circ \hat f_{i_1 i_2} |_{\DD_{i_1 \ldots 
i_n}}
\, .
$$

The map $\hat f^{(n)}_{\vec \imath}$ has a fixed 
point in $\DD_{i_1 \ldots i_n}$
if and only if $i_{n+1}=i_1$, the fixed
point is then unique. It follows that
the ``hyperbolic points of period $n$ for $\hat f$,'' i.e., those 
$(w_1, w_2)\in \DD_{\vec \imath}$, for some $\vec \imath\in \SS^{n+1}\setminus\SS_0^{n+1}$
such that
$f^{(n)}_{\vec \imath} (w_1, w_2) =(w_1, w_2)$, are in bijection with
the  hyperbolic symbolic ``cycles'' of length $n$, i.e., admissible sequences
$\vec \jmath \in \SS^n\setminus \SS_0^{n}$ such that $t_{j_n j_1}=1$.
\endproclaim

\subhead 2.B Banach spaces and elementary transfer operators
\endsubhead

{\bf The Banach spaces $\BB_k$, $\BB_k'$}

We now associate  Banach spaces of complex functions to
each $k \in \SS_0 \cup S_1$. First, for each $k\in  \SS_1$
we set  $\BB'_k=\BB_k$ to be the Banach space
of  holomorphic functions on the interior of $(\overline \complex \setminus 
\DD^1_k) \times \DD^2_k$, vanishing at $\{\infty \}\times \DD^2_k$, 
and  which extend
continuously to the boundary, with the supremum norm.

For $k \in \SS_0$, we distinguish between cases (P.a) and (P.b).
In  case (P.a), we shall define below an open simply
connected subset $U^2_k$ of $\DD^2_k$,
containing the compact set
$$
K^2_k :=  { \cup_{j\in \SS_1}  \phi_{kj,s} 
(\DD^1_k, \DD^2_j) } \, , \tag{2.11}
$$
and such that $\phi_{kk,s}(w_1, U^2_k)
\subset U^2_k$ for all $w_1\in  \partial \DD^1_k$. 
The space $\BB_k$ will be a
subset of the space of analytic functions in  $(\bar \complex
\setminus\DD^1_k) \times U^2_k$, vanishing at $\{\infty\} \times U^2_k$, 
and  extending continuously
to $\partial \DD^1_k \times U^2_k$, endowed with a norm to
be introduced below, using Fatou coordinates. We shall use in the proof
of Lemma~3.8 that the Banach
space $\BB'_k$
of holomorphic functions in  $(\bar \complex \setminus\DD^1_k) \times
U^2_k$, which vanish at $\{\infty\} \times U^2_k$ and
extend continuously to the
boundary, is continuously embedded in $\BB_k$. 

In case (P.b),
we shall introduce below an  open simply connected
subset  $U^1_k $ of the attracting petal  $-\DD^1_k$  which does not intersect
the compact set
$$
G^1_ k:=
\cup_{j\in \SS_1} (\phi_{jk, u} (\DD^1_j, \DD^2_k))  \, ,\tag{2.12}
$$
and such that $\phi_{kk,u}^{-1}(U^1_k,z_2) \subset U^1_k$ for all $z_2 \in \DD^2_k$
(recall Definition ~2.5).
The Banach space $\BB_k$ 
will be a subset  of the analytic functions in  the interior of $U^1_k \times
\DD^2_k$, extending continuously to $U^1_k \times
\partial \DD^2_k$, endowed with a norm to be defined below.
We shall use in  Lemma~3.8 that the Banach space $\BB'_k$ of holomorphic
functions in  $U^1_k \times \Int \DD^2_k$
which extend continuously to the boundary is continuously embedded
in $\BB_k$. Note that  case (P.b) did not occur
in the one-dimensional situation studied by Rugh  [Ru3].

\medskip
{\bf The elementary transfer operators $\widehat \LL_{kj}$}

Next, to each  $j, k$ with $t_{kj}=1$,
we associate an {\it elementary transfer operator}.
If $(k ,j)\in \SS_1\times  (\SS_1\cup \SS_0)$, we set for $\psi \in \BB_k$:
$$
\eqalign{
&\widehat \LL_{kj} \psi \, (z_1, z_2) =
\oint_{\partial  \DD^1_k} \oint_{\partial   \DD^2_k}
{dw_1 \over 2i\pi} {dw_2 \over 2i\pi}  
{s_{\phi'_{kj,s}} \partial_2 \phi_{kj,s} (w_1, z_2) 
\over
w_2 - \phi_{kj,s} (w_1, z_2)}
{\psi(w_1,w_2) \over z_1 - \phi_{kj,u} (w_1, z_2)}  \, .
\cr&
}\tag{2.13}
$$
In \thetag{2.13}, and from now on, $s_{\phi'_{kj,s}}$ is the (well-defined) sign of 
of $\partial_2 \phi_{kj,s}$ on $\II^1_k\times \II^2_{j}$.
See Remark~3.6 below for a heuristic justification of this choice, and pages 1248--1250 in [Ru1]
for a more analytic explanation.

\smallskip

If $(k ,j)\in \SS_0\times  \SS_1$, we use \thetag{2.13},
replacing $\partial  \DD^2_k$ by a simple curve $\Gamma^2_k$ inside 
$U^2_k$ which does not intersect
$K_k^2$  in case (P.a) (this ensures that $\psi$ is
well defined on $(\partial \DD^1_k,\partial  \Gamma^2_k)$ and the
integral is holomorphic).
If $k$ is of type (P.b), we replace
$\partial  \DD^1_k$ by a curve $\Gamma^1_k$ inside $U^1_k$ 
such that (see Definition~2.5)
$$
\cup_{j \in \SS_1}(\phi_{ k j, u} (\Gamma^1_j, \DD^2_k))
\cap -\DD^1_j =\emptyset\, .
$$ 

\smallskip

If $k\in \SS_0$ is of type  (P.a),  we define $\widehat \LL_{kk} \psi \, (z_1, z_2)$
for $(z_1, z_2)\in (\bar \complex \setminus\DD^1_k) \times U^2_k$
by
$$
\eqalign{
&\widehat \LL_{kk} \psi \, (z_1, z_2) =
\oint_{\partial \DD^1_k}  {dw_1 \over 2i\pi}   
s_{\phi'_{kk,s}} \partial_2 \phi_{kk,s} (w_1, z_2) 
{\psi(w_1,\phi_{kk,s} (w_1, z_2)) \over z_1 - \phi_{kk,u} (w_1, z_2)} \,  .
}\tag{2.14}
$$

In case (P.b), 
for $(z_1,z_2)\in U^1_k \times \DD^2_k$,  we use 
$$
\eqalign{
\widehat \LL_{kk} \psi \, (z_1, z_2) &=
s_{\phi'_{kk,s}}  \cdot \cr
&{\partial_2 \phi_{kk,s} (\phi_{kk,u}^{-1}(z_1, z_2),z_2) 
\over  
 \partial_1 \phi_{kk,u}(\phi_{kk,u}^{-1}(z_1, z_2), z_2)}
\psi(\phi_{kk,u}^{-1}(z_1, z_2),\phi_{kk,s} (\phi_{kk,u}^{-1}(z_1, z_2), 
z_2)) \, .  }\tag{2.15}
$$
Note that 
$$
\widehat \LL_{kk} \psi \, (z_1, z_2)=
\oint_{\partial   \DD^2_k} {dw_2 \over 2i\pi}  
{s_{\phi'_{kk,s}} \partial_2 \phi_{kk,s} (\phi_{kk,u}^{-1}(z_1, z_2), z_2) 
\over \partial_1 \phi_{kk,u}(\phi_{kk,u}^{-1}(z_1, z_2), z_2)}
{\psi(\phi_{kk,u}^{-1}(z_1, z_2),w_2)  \over
w_2 - \phi_{kk,s} (\phi_{kk,u}^{-1}(z_1, z_2), z_2)} 
\, .\tag{2.16}
$$
(The two signs cancels in the
residue computation because $z_1$ is outside of $\DD^1_k$.)

\remark{Remark 2.9}
\roster
\item
If  $j\in \SS_1$, then $\widehat \LL_{kj}$ is bounded from $\BB_k$ to $\BB_j=\BB_j'$. 
\item
If $k\in \SS_1$ and $j \in \SS_0$ then $\widehat \LL_{kj}$ is bounded 
from $\BB_k$ to $\BB'_j$, and $\BB'_j \subset \BB_j$, continuously
by our
conditions on $U^1_j$, $U^2_j$ in cases (P.a), respectively (P.b),
so that $\widehat \LL_{kj}$ is bounded from $\BB_k$ to
$\BB_j$. 
\endroster
\endremark

\smallskip

{\bf Approximate Fatou coordinates, more about Banach spaces}

Our next step is to give a precise definition of the
Banach spaces $\BB_k$ associated to $k\in \SS_0$ and to study the corresponding
elementary transfer operators.

Let $\hat f_{kk}$ be  of type (P.a).
It is well-known 
that the injective map $\FF(z)=1/z$  is an approximate Fatou coordinate, i.e.,
$$
\eqalign
{
&\FF( \phi_{kk,s}(w_1,z_2)) =\FF(z_2)+1+
z_2 \cdot \EE_{k,2}(w_1,z_2)  \, , \cr
&\EE_{k,2} : \widetilde \DD^1_k\times \DD^2_k \to \complex\, ,  {\text{ holomorphic and bounded.} 
}
}\tag{2.17}
$$
The set $\Omega^2_k=\FF(\Int \DD^2_k)\subset \complex$ is open and simply connected.
It is easy to check that $\Omega^2_k+1$ is contained
in $\Omega^2_k$ and that for some (large) $R_k$ the domain
$\Omega^2_k$ contains the closed ``right''  half plane $\overline{ H_{R_k}}$ where
$$H_R=\{z \in \complex \mid \Re z > R\} \, .$$

Note that the   one-dimensional Fatou coordinates
$\FF_{w_1}$ associated to
each $\phi_{kk,s}(w_1,\cdot)$ (see e.g [Mi] or [Ru3, Lemma 2.1])
solves
 $$\FF_{w_1}( \phi_{kk,s}(w_1,z_2)) =\FF_{w_1} (z_2)+1 \, ,\tag{2.18}
$$
while we would ``like'' 
$\widetilde \FF_{w_1}( \phi_{kk,s}(w_1,z_2)) =\widetilde \FF_{\phi_{kk,u}(w_1,z_2)} (z_2)+1 $,
which is not immediately available in the attracting petal. 
Our argument will be perturbative --- in Lemma~ 3.2 we shall compare our elementary
operator $\widehat \LL_{kk}$ to
a direct product --- it hence is possible to work with the approximate
Fatou coordinate $\FF$, which also has the important feature of
being holomorphic on $\complex^*$.

\smallskip

In  case (P.b), recall from Definition~2.5 that the map $\phi_{kk, u}^{-1} (\cdot, z_2)$ is the inverse of
$w_1\mapsto \phi_{kk, u} (w_1, z_2)$ for fixed $z_2\in \widetilde \DD^2_k$
and $w_1 \in -\DD^1_k$.
 We have
$$
\eqalign
{
&\FF(\phi_{kk,u}^{-1} (z_1,z_2)) =\FF(z_1)-1+  z_1\cdot\EE_{k,1}(z_1,z_2)\,
, \cr
&   \EE_{k,1} : - \DD^1_k\times \widetilde \DD^2_k \to \complex \, , 
{\text{ holomorphic and bounded}}  \, . }\tag{2.19}
$$
Here  $\Omega^{1-}_k=\FF(\Int(-\DD^1_k))$ is open, simply connected, and  there is $R_k< 0$  so that
$\Omega^{1-}_k$ contains the closed left half plane $-\overline{ H_{-R_k}}$.

\definition{Definition of the Banach space $X(H_R)$}
For $R > 0$, let $X(H_R)$ be the isometric image of $L^1(\real_+, \hbox{Lebesgue})$ 
under the shifted Laplace transform
$$
\tilde \psi(\tilde w)= 
\int_0^\infty \psi_L(t) e^{-(\tilde w-R)t} \,  dt \, ,
\, \tilde w \in H_R \, , \tag{2.20}
$$
inside 
the space of holomorphic functions in $H_R$, with induced norm.
\enddefinition 

(We refer to Doetsch [Doe] for the basics of the Laplace transform.)
Functions in $X(H_R)$ are in fact bounded in $H_R$ by the $L^1$ norm
of $\psi$. One can easily check (see e.g. Lemma 2.5 in [Ru3] for similar
ideas) that for any closed right-half plane
$H_{R'}$ with $R' > R$ there is a constant $C_{R'-R}$ so that
for each $\tilde \psi \in X(H_{R})$ the derivative $\tilde \psi'$ is bounded
in $H_{R'}$ by $C \cdot \int |\psi_L(t)| \, dt$, where $\psi_L$ is
an $L^1$ representant of the inverse Laplace transform
of $\tilde \psi$ in $H_R$ (use that $e^{-\delta t} t$ is bounded
on $[0, \infty]$). It is not
difficult  [Ru3, Lemmas 2.2 and 2.3] to prove that the spectrum of the
translation operator $S: X(H_R )\to X(H_R)$ defined by $S\tilde \psi(w)
=\tilde \psi(w+1)$ is the unit interval $[0,1]$. 

\smallskip
 
{\bf Laplace coordinates and $\BB_k$ norm in case (P.a)}

We next exploit the Fatou coordinates, adapting some  
definitions from [Ru3]. Let us consider first
the case (P.a). Although $\Omega^2_k$ contains the
closed right half-plane $H_{R_k}$, we shall need to work with a
slightly larger domain.   The set
$K^2_k $ from \thetag{2.11}
is a compact subset of  $\Int  \DD^2_k$.
We may thus adapt  Lemma 2.4 in [Ru3], together with
the arguments presented just after it,   finding
$R_k > m_k>0$ and an open connected and simply connected subset
$N^2_k$ of $\Omega^2_k$
so that
$$
\eqalign{
&\overline {H_{ R_k} }\subset N^{2}_k\subset H_{R_k} - m_k 
\, , \quad K^2_k \subset \FF^{-1}  (N^2_k) \, , \cr
& N^{2}_k+1 \pm {\sup|\EE_{k,2}| \over |z|}\subset N^{2}_k  \, , \forall z \in H_{R_k-m_k}  \, .\cr }
$$

\definition{Definition ($X(N^2_k)$)}
For $\widetilde R_k > R_k$, we let $X(N^2_k)$ be the subset of $X(H_{\widetilde R_k})$
consisting in those functions which admit an analytic continuation
to $N^2_k$ with a continuous extension to the boundary.
We take as norm the sum of the supremum norm on $N^2_k$ 
with the $X(H_{\widetilde R_k})$ norm of the restriction to
$H_{\widetilde R_k}$.
\enddefinition 

Strictly speaking, we should replace our translation
operator $S$  on $X(H_{\widetilde R_k})$ by a translation operator $T$ which
``lives'' in  $X(N^2_k)$. Since this does not
influence the spectrum (details are to be found
in [Ru3, (2.27--2.30)]) we shall instead
abuse notation.

\definition{Definition of $U^2_k$ and $\BB_k$, $\widetilde \BB_k$}
Define an open subset of $\DD^2_k$ by 
$
U^2_k =\Int  \FF^{-1}_2 (N^2_k) 
$, it has the required properties.
We use the notation $\AA(\bar \complex \setminus \DD^1_k)$ for the Banach
space of holomorphic functions in  $\bar \complex \setminus \DD^1_k$
vanishing at infinity and extending continuously to the boundary.
Let $\BB_k$ be the set of analytic functions
in $(\bar \complex \setminus \DD^1_k)\times U^2_k$  
given by the isometric image of the Banach space tensor product 
$
\widetilde \BB_k=\AA(\bar \complex \setminus \DD^1_k) \otimes  X(N^2_k)
$
under $J^{-1}$ defined by
$$
\psi= ( J^{-1} \tilde \psi) \, , \quad
\psi(z_1, z_2)= -\tilde \psi(z_1, \FF(z_2))
\cdot  \FF'(z_2)=\tilde \psi(z_1, \FF(z_2)) z_2^{-2}\, , \tag{2.21}
$$
with induced norm.  In particular,
if $\tilde \psi \in \widetilde \BB_k$ then
$\tilde \psi (w_1, \cdot) \in X(N^2_k)$ for every $w_1 \in \bar \complex \setminus \DD^1_k$.
\enddefinition

Writing $\EE'_{k,2}= \partial_2 \EE_{k,2}$,
the elementary transfer operator $\widehat \LL_{kk}$ \thetag{2.14}
in the Fatou coordinates
(i.e.,  acting on  $\widetilde \BB_k$
via \thetag{2.21}) can be written as
$$
\eqalign
{
&\widetilde \LL_{kk} \tilde \psi(z_1,\tilde z_2)
= 
s_{\phi'_{kk,s}} \oint_{\partial D^1_k}\cr
&\,\,  {dw_1\over 2i\pi}
\biggl (1- {\EE_{k,2}(w_1,1/\tilde z_2)+\tilde z_2^{-1} \EE_{k,2}'(w_1,1/\tilde z_2)
\over \tilde z_2 ^2} \biggr ) \cdot
{\tilde \psi(w_1,\tilde z_2+1+\tilde z_2^{-1}\EE_{k,2}(w_1,1/\tilde z_2)) 
\over z_1 - \phi_{kk,u} (w_1, \FF^{-1}(\tilde z_2))} \, .} \tag{2.22}
$$
Indeed, differentiating both sides of \thetag{2.17} with respect
to $z_2$ yields 
$$
\FF'(\phi_{kk,s}(w_1,z_2)) 
\partial _2\phi_{kk,s} (w_1, z_2) = 
\FF'(z_2) \biggl (1+ { \EE_{k,2}(w_1,z_2)+z_2 \EE_{k,2}'(w_1,z_2)\over \FF'(z_2) } \biggr ) \, , 
$$
so that $\widehat \LL_{kk} ( J^{-1}\tilde \psi) $ coincides with $J^{-1} 
(\widetilde \LL_{kk} \tilde \psi)$, using \thetag{2.17} again.

\smallskip

{\bf Laplace coordinates and $\BB_k$ norm in case (P.b)}

Let us now discuss  case (P.b). We shall use  left half-planes 
$-H_{-R}$ for $R< 0$
and spaces  $X(-H_{-R})$ of ``left'' Laplace transforms
$$
\tilde \psi(\tilde z)= 
\int_0^\infty \psi_L(t) e^{(\tilde z-R)t} \,  dt \, , \, 
\Re \tilde z  < R \, , \, \psi_L\in L^1(\real^+, \, \hbox{Lebesgue}\, ) \, .   \tag{2.23}
$$

Since the closed attracting petal $-\DD^1_k$ does not intersect 
$G^1_k$ from \thetag{2.12},
the open left half plane $N^1_k=-H_{-R_k}$ satisfies
$\FF^{-1}(N^1_k)  \cap G^1_k =\emptyset$ and
$$
N^1_k -1 \pm {\sup |\EE_{k,1}|\over |z|}\subset N^1_k \, , \forall z \in -H_{-R_k} \, .
$$

\definition{Definition of $X(N^1_k)$, $U_k^1$, $\BB_k$, $\widetilde \BB_k$}
Fix $\widetilde R_k <  R_k < 0$. Let $X(N^1_k)$ be the subset of  $X(-H_{-\widetilde R_k})$ 
consisting in those functions which extend analytically
to $N^1_k$ and continuously to its boundary, with norm 
the sum of
the $L^1$ norm of the inverse Laplace transform
with the supremum  in $N^1_k$.
Set
$
\widetilde \BB_k= X(N^1_k) \otimes \AA(\DD^2_k)
$
(with $\AA(\DD^2_k)$ the Banach space of holomorphic functions
in $\Int \DD^2_k$ extending continuously to the boundary).

Defining the following analogue of \thetag{2.21}
$$
\eqalign
{
\psi&= ( J^{-1}\tilde \psi) \, , \quad
\psi(z_1, z_2)= -\tilde \psi(\FF(z_1),z_2) \FF'( z_1) \, , \cr
}\tag{2.24} 
$$
we finally let
$\BB_k$ be the isometric image of
$\widetilde \BB_k$ under $J^{-1}$.  
\enddefinition

Thus, elements of $\BB_k$ are analytic
functions in the interior of $U^1_k\times \DD^2_k$
for the open connected set
$U^1_k=  \FF^{-1}(N^1_k)  $.

Differentiating \thetag{2.19} with respect to $z_1$
(writing $\EE'_{k,1}=\partial_1 \EE_{k,1}$) yields
$$
\FF'(\phi_{kk, u}^{-1}(z_1,z_2)) \partial_1\phi_{kk, u}^{-1}(z_1,z_2)=
\FF'(z_1) \biggl (1+ {\EE_{k,1}(z_1,z_2) + z_1 \EE_{k,1}'(z_1, z_2)\over
\FF'(z_1)}\biggr ) \, .
$$
Set $\upsilon_1= \FF^{-1}(\tilde z_1 -1 +\tilde z_1^{-1} \EE_{k,1}(\tilde 
z_1^{-1},z_2) )$ and  note that 
$$
\upsilon_1=\phi_{kk,u}^{-1}(\FF^{-1} \tilde z_1,z_2)
\, . \tag{2.25}
$$
The transfer
operator  \thetag{2.15} in the Fatou coordinates can be written for 
$(\tilde z_1, z_2) \in  \FF (U^1_k) \times \DD^2_k$ as
$$
\eqalign{
\widetilde \LL_{kk} \tilde \psi(\tilde z_1, z_2)
&= s_{\phi'_{kk,s}} \biggl (1 -
{\EE_{k,1}(1/\tilde z_1,z_2)+
\tilde z_1^{-1} \EE_{k,1}' (1/\tilde z_1,z_2)\over \tilde z_1 ^2} \biggr )\cr
&\qquad  \qquad\qquad \quad\oint_{\partial   \DD^2_k} {dw_2 \over 2i\pi}  
{\partial_2 \phi_{kk,s} (\upsilon_1, z_2) 
\over
w_2 - \phi_{kk,s} (\upsilon_1, z_2)}
\tilde \psi(\FF \upsilon_1 ,w_2)
\, . }\tag{2.26}  
$$
(Use
$( \partial_1 \phi_{kk,u}(\phi_{kk,u}^{-1}(z_1, z_2), z_2))^{-1}
=  \partial_1 \phi_{kk,u}^{-1}(z_1, z_2)$.)

\head 3. Spectrum and determinants  of the symbolic transfer operator
\endhead

\subhead 3.A  Direct product transfer operators for (P.a)--(P.b)
\endsubhead

For $k\in \SS_0$ of type (P.a), we introduce a
direct tensor product operator, written using the Laplace transforms 
\thetag{2.20} of $\tilde \psi(w_1, \tilde w_2)$ for each fixed $w_1$ as 
$$
\widetilde \LL_{kk}^\otimes \tilde \psi(z_1, \tilde z_2)=
s_{\phi'_{kk,s}}
 \oint_{\partial D^1_k}\,  {dw_1\over 2i\pi}
{1\over  ( z_1 - \lambda_{k,u} w_1)}
\int_0^\infty  e^{-(\tilde z_2-\widetilde R_k +1)t} \, \psi_L(w_1,t)  \,  dt \, ,\tag{3.1}
$$
where $\lambda_{k,u} = \partial_1 \phi_{kk,u} (0, 0)$.

If $k$ is of type (P.b), setting $\lambda_{k,s}={\partial_2 \phi_{kk,s} (0, 0)}$,
the corresponding direct tensor product approximation can be written as
$$
\widetilde \LL_{kk}^\otimes \tilde \psi(\tilde z_1, z_2)
= s_{\phi'_{kk,s}} 
{\lambda_{k,s} }
\int_0^\infty  e^{(\tilde z_1-\widetilde R_k -1)t} \psi_L(t,\lambda_{k,s} z_2)  \, dt 
 \,  .\tag{3.2}
$$

We define a direct sum of direct
tensor products  $\widetilde \LL_0^\otimes$ acting on
$\widetilde \BB_0=\oplus_{k\in \SS_0}  \widetilde \BB_k$ by setting
$\widetilde \LL_0^\otimes=\oplus_{k\in \SS_0} \widetilde \LL_{kk}^\otimes$.
We use similar notations $\BB_0$,
$\LL_0^\otimes$, corresponding to the conjugated  operators $\LL_{kk}^\otimes$ on $\BB_k$.

\proclaim{Lemma 3.1
(Spectrum and resolvent of the direct products $\widetilde \LL_0^{\otimes}$,
$\LL_0^{\otimes}$)}
\roster
\item
Let $\{ \lambda_{k, u}\}$, in case (P.a), and $\{ \lambda_{k, s}\}$ 
in case (P.b), be defined above.
The operator $\widetilde \LL_0^\otimes$ is bounded on $\widetilde \BB_0$
and its spectrum is 
is the following set:
$$
\eqalign{
&\{ [0,1] \, , k \hbox{ type (P.a)}\, , \lambda_{k, u} > 0\} \cup
\{[\lambda_{k, u} ,1] \, , k \hbox{ type (P.a)}\, , \lambda_{k, u} < 0   \}\cr
& \cup
\{[0,  \lambda_{k, s}] \, ,  k \hbox{ type (P.b)} \, , \lambda_{k,s} > 0  \}\cup
\{[-\lambda_{k,s}^2, -\lambda_{k,s} ] \, , k \hbox{ type (P.b)} \, , \lambda_{k,s} < 0  \}
\, . }
$$
\item
Let $\widetilde \BB_0(\epsilon)$ be the Banach space obtained by replacing discs
$\DD^{1,2}_k$ of radius $r$
by discs of radius $\epsilon r$, and,
in case (P.a) $\widetilde R_k$ by $\widetilde R_k+ 1/\epsilon$, 
$N^2_k$ by $N^2_k + 1/\epsilon$, in case (P.b) $\widetilde R_k$ by $\widetilde R_k- 1/\epsilon$,
$N^1_k$ by $N^1_k -1/\epsilon$.
Then, for every $1/z\notin \sp\widetilde \LL_0^\otimes$ (on $\widetilde \BB_0$) there is
$C(z)$ so that
$$
\|(1-z \widetilde \LL_0^\otimes)^{-1}\|_{\widetilde \BB_0(\epsilon)}
\le C(z)  \|(1-z\widetilde \LL_0^\otimes)^{-1}\|_{\widetilde\BB_0} \, , \quad \forall \epsilon > 0\, .
$$
\endroster
\therosteritem{1--2} hold for $\LL_0^\otimes$ acting on $\BB_0$ and $\BB_0(\epsilon)$.
\endproclaim

Note that in the definition of $\BB_0(\epsilon)$ we  disregard the condition
$K^2_k \subset \FF^{-1}(N^2_k)$, but we do use, in case (P.a) e.g., that condition (H.2) also holds
for the disc $\epsilon \cdot \DD^1_k$ and the restricted domain $\DD^2_k
\cap \FF^{-1}(N^2_k +1/\epsilon)$. See also Remark~4.7.

\demo{Proof of Lemma 3.1}
Since
$\widetilde \LL_0^{\otimes} $ is a direct sum, it
suffices to consider each term $\widetilde \LL_{kk}^\otimes$ 
acting on $\widetilde \BB_k$. So, let us fix
$k\in \SS_0$, assuming first that we are in case (P.a).
Then,  since we assumed that the neutral eigenvalue is
$+1$, we have $s_{\phi'_{kk,s}}=1$.
The operator
$\widetilde \LL^{\otimes}_{kk}$  \thetag{3.1} is thus the direct tensor product
$
\widetilde \LL^{\otimes}_{kk}= \MM_{k,u}  \otimes \TT^+ 
$,
with
$$
\eqalign{
&\MM_{k,u} \psi_1(z_1)=
\oint_{\partial D^1_k}\,  {dw_1\over 2i\pi}
{1\over  ( z_1 - \lambda_{k,u} w_1)}\psi_1(w_1)  \, , \cr
&\TT^+\widetilde \psi_2 (\tilde z_2) = \int_{0}^\infty 
e^{-(\tilde z_2-\widetilde R_k+1)t} \psi_L(t) \, dt \, . }
$$
The results of Ichiniose [Ic]
give that its spectrum is just $\{ \sigma_1 \cdot \sigma_2 \, ,
\sigma_1 \in \sp(\MM_{k,u})\, ,
\sigma_2 \in \sp(\TT^+)\}$. The remarks above and the ideas
in [Ru3] easily yield that the spectrum of $\TT^+$ on $X(N^2_k)$ is
the unit interval $[0,1]$. One also obtains that
$\MM_{k,u}$ is nuclear on the Banach space $\AA(\bar \complex \setminus \DD^1_k)$ 
of functions holomorphic in $\bar \complex \setminus \DD^1_k$,
vanishing at $\infty$, and extending continuously to $\partial \DD^1_k$, and 
that its spectrum consist of the  eigenvalues
$\lambda_{k,u}^\ell$ for $\ell \in \integer_+$. (Indeed, the trace of
$\MM_{k,u}^\ell$ can be shown, following the methods
in [Ru1, Ru2, Ru3] to be
$1/(1-\lambda_{k,u}^\ell)$.)

If we are in case (P.b), we proceed similarly, starting
from \thetag{3.2} and using
$$
\eqalign
{&\TT^-\widetilde \psi (\tilde z_1) = \int_{0}^\infty 
e^{(\tilde z_1-\widetilde R_k-1)t} \psi_L(t) \, dt \, , 
\cr
&\MM_{k,s} \psi_2(z_2)
= \oint_{\partial D^2_k}\,  {dw_2\over 2i \pi}
{s_{\phi'_{kk,s}} \cdot \lambda_{k,s}  \over  (w_2 - \lambda_{k,s} z_2)}
\psi_2(w_2)  = | \lambda_{k,s} | \cdot \psi_2(\lambda_{k,s} z_2)  \, . }
$$
For this, just check that the trace of
$\MM_{k,s}^\ell$ is
$|\lambda_{k,s}^{\ell}|/(1-\lambda_{k,s}^{\ell})$. Its
spectrum is thus  
$(\hbox{sgn} \, \lambda_{k,s}) \cdot (\lambda_{k, s})^{\ell}$ for  $\ell \in \integer^*_+ $.

\smallskip

For the second claim, 
observe first that modifying $\widetilde R_k$ and $N^i_k$
does not affect the spectrum or the norm of the resolvent,
by definition of the Laplace norm and because
$\FF'(w_1)/\FF'(w_1\pm 1/\epsilon)$ is uniformly bounded.
It thus suffices to consider $\widetilde \LL_0^\otimes$
on  $\widetilde \BB_0(\epsilon)$, where only the
hyperbolic radii of $\DD^{1,2}_k$ have been rescaled.

Then,  denoting by $\HH_\epsilon: \widetilde \BB_0(\epsilon)
\to \widetilde \BB_0$ the following isometry  on
$\widetilde \BB_k$ for $k\in \SS_0$: $\HH_{\epsilon, k} \varphi(w_1, \tilde w_2)
=\varphi (\epsilon w_1, \tilde w_2)$, and similarly in case (P.b), we see  that
$\HH_\epsilon^{-1} \widetilde \LL_0^\otimes \HH_\epsilon=\widetilde \LL_0^\otimes$
(this follows from the facts that the ``weight'' in $\MM_{k, (u,s)}$ is constant
while the ``dynamics'' is linear) and thus
$$
\HH_\epsilon^{-1} (1-z \widetilde \LL_0^\otimes) ^{-1}\HH_\epsilon
= (1-z \widetilde \LL_0^\otimes) ^{-1} \, .
$$
Since $\|\HH^{\pm 1}_\epsilon \|=1$ for all $\epsilon$, we have proved the bounds
on $\widetilde \BB_0(\epsilon)$.
\qed
\enddemo

\proclaim{Lemma 3.2 (Perturbation theory for $\widehat \LL_0$)}
$\widehat \LL_0$ is bounded on $\BB_0$ and on $\BB_0(\epsilon)$, and
there is a positive function $\GG$ with $\GG(\epsilon)\to 0$
as $\epsilon \to 0$ such that
$$
\sp (\widehat \LL_0|_{\BB_0(\epsilon)}) \subset \{ z \in \complex \mid
\exists \tilde z \in \sp \widetilde \LL_0^\otimes \, , | z - \tilde z | =\GG(\epsilon) \}
 \, . 
$$ 
\endproclaim

We summarize in two sublemmas  elementary properties of
the Laplace  norm which will be used in the proof
of Lemma~3.2:

\proclaim{Sublemma 3.3}
Let $\tilde \psi_1 \in X(H_R)$ for $R > 0$.
 Then:
\roster
\item
For  each integer
$s \ge 1$, the function $\tilde \psi_1(\tilde z)/\tilde z^s$ belongs to
$X(H_R)$ and 
$$\|\tilde \psi_1(\tilde z)/\tilde z^s\|_{X(H_R)} \le
\|\tilde \psi_1\|_{X(H_R)}/R^s\, .$$

\item
For every $w$ such that $\Re w > 0$, the functions $\tilde \psi_1(\tilde z+w)$
and $\tilde \psi'_1(\tilde z+w)$
belong to $X(H_R)$ and 
$$\max \bigl ( \|\tilde \psi_1(\tilde z+w)\|_{X(H_R)},\|\tilde \psi'_1(\tilde z+w)\|_{X(H_R)}\bigr )
\le \|\tilde \psi_1\|_{X(H_R)} /(e \cdot \Re w)\, .
$$

\item Assume that $\tilde \psi_1$ admits a bounded holomorphic
extension to $H_{R-2\eta}$ for $\eta > 0$ and  $R-\eta> 0$.
Let  $\Phi : H_{R-2\eta} \to H_{R-2\eta}$ 
be holomorphic. For  each integer
$s \ge 2$, the function $\tilde \psi_1(\Phi(\tilde z))/\tilde z^s$ belongs to
$X(H_R)$ and 
$$
\biggl \| {\tilde \psi_1(\Phi(\tilde z))\over \tilde z^s} \biggr \| _{X(H_R)}\le
\sup_{H_{R-2\eta}} |\tilde \psi_1|{4 \over  \eta \cdot (R-\eta)^{s-1}} \, .
$$

\item
If $\EE(z)$ is holomorphic and bounded in a neighbourhood  $\DD$ of zero
such that $\FF(\DD)$ contains a
closed half-plane $H_{R'}$, $0<R' < R$,
then $\EE(1/\tilde z) \tilde \psi_1(\tilde z) \in X(H_R)$ and 
$$
\|\EE(1/\tilde z) \tilde \psi_1(\tilde z) \|_{X(H_R)}
\le {2\pi \over 1-R'/R} \sup_{\DD } |\EE|  \cdot   \| \tilde \psi_1\|_{X(H_R)} \, .
$$

\endroster
(There are obvious analogues for $\tilde \psi_2 \in X(-H^-_{R})$ with $R < 0$.)
\endproclaim 

\proclaim{Sublemma 3.4}
Assume that $\tilde \psi_1 \in X(H_R)$ for $R > 0$, and
that $\tilde \psi_1$ admits a bounded holomorphic
extension to $H_{R-2\eta}$ for some $\eta > 0$. 
If $\EE(z_1,z_2)$ is analytic and bounded in 
$\DD^1 \times \DD^2$ such that $\FF(\DD^1)$ 
contains a closed half-plane $H_{R'}$ with $0<R' < R$,
and $\DD_2$ is a neighbourhood of zero,
then for each $\rho <1$
the product  $\EE(1/\tilde z_1,z_2) \cdot \tilde \psi_1(\tilde z_1)$ is an 
element of 
$X(H_R)\otimes \AA(\rho \DD^2)$ and 
$$
\|\EE(1/\tilde z_1,z_2) \cdot  \tilde \psi_1(\tilde z_1) \|
\le \sup_{\DD^1\times \DD^2 }|\EE| { \| \tilde \psi_1\|_{X(H_R)}
\over (1-\rho) (1-R'/R) }\, .
$$
\endproclaim

\demo{Proof of Sublemma 3.3}
(1)
By induction, it is enough to consider  $s=1$.
We first note that $\tilde z^{-1} =
\int_0^\infty e^{-(\tilde z-R)t } e^{-Rt}\, dt$ belongs to
$X(H_R)$. Since $e^{-Rt}$ is also  in $L^\infty$,
we may obtain an ($L^1$) inverse Laplace transform
of $\tilde \psi_1(\tilde z)/\tilde z$ by performing
the convolution
$$
e^{-Rt} \star \psi_{1,L}(t) =
\int_0^t e^{-R\tau} \psi_{1,L}(t-\tau)d\tau \, .
$$
Now, 
$
\int_0^\infty |e^{-Rt} \star \psi_{1,L}(t)| dt
\le\int_0^\infty d\tau e^{-R\tau} 
\int_\tau^\infty dt |\psi_{1,L}(t-\tau)|
\le \|\tilde \psi_1\| \int_0^\infty e^{-R \tau}d\tau $
which is equal to
$\|\tilde \psi_1\|_{X(H_R)}/R$.

(2) By definition
$
\tilde \psi'_1(\tilde z+w)=
-\int_0^\infty t e^{-(\tilde z + w -R)t}\psi_L(t) \, dt
$.
Estimating the $L^1$ norm of
$-t e^{-wt} \psi_{1,L}(t)$, we obtain the claim for $\tilde \psi_1'(\tilde z+w)$ the
bound for $\tilde \psi_1(\tilde z+w)$ is easier. 

(3) We can express the
inverse Laplace transform of $\tilde \psi_1(\Phi(\tilde z))/\tilde z^s$ as
$$
\eqalign
{
\psi_{1,L,s}(t)&=
\int_{\partial H_{R-\eta}} {du\over 2i\pi} e^{t(u-R)}{\tilde \psi_1(\Phi(u))\over  u^s}
= e^{-t\eta}\int_{\partial H_{R}} {dv\over 2i\pi} e^{t(v-R)}
{\tilde \psi_1(\Phi(v-\eta))\over (v-\eta)^s}\, .
}
$$
It follows that 
$$
\eqalign
{
|\psi_{1,L,s}(t)| \le e^{-\eta t} \sup_{H_{R-\eta}}|\tilde \psi_1| 
\int_{-\infty}^\infty{dx \over 2\pi} {1\over ((R-\eta)^2+x^2)^{s/2}}\, ,
}
$$
which gives the claimed bound.

(4) 
Let $\EE(z)=\sum_{m=0}^\infty a_m z^m$ be a Taylor series
for $\EE$ at the origin. By the Cauchy formula we have
 $| a_m|\le 2\pi \sup_{\DD} |\EE| (R')^m$ for $m\ge 0$.
Then
$$
\EE(1/\tilde z) \tilde \psi_1(\tilde z)
=
\sum_{m= 0}^\infty a_m {\tilde \psi_1(\tilde z)\over \tilde z^m}\, .
$$
Then, bound \therosteritem{1}  shows that the sum is 
$\le \sum_{m\ge 0} |a_m|\|\tilde \psi_1\|  R^{-m}$.\qed
\enddemo

\demo{Proof of Sublemma 3.4}
As in the proof of Sublemma 3.3  \therosteritem{4}, we  use
a Taylor series
$$
\EE(z_1, z_2)=\sum_{m=0}^\infty a_m(z_1) z_2^m \, ,
$$
where each $a_m(z_1)$ is holomorphic in $\DD^1$
and, if $\delta>0$ is the radius of $\DD^2$, 
$$
 \sup_{\DD^1} |a_m| \le 
{2 \pi \sup_{\DD^1\times \DD^2} | \EE|
\over  \delta^m } \, .
$$
By   Sublemma 3.3 \therosteritem{4}, we know that
each $a_m(1/\tilde z_1) \tilde \psi_1(\tilde z_1)$ belongs
to $X(H_R)$ with 
$\|a_m(1/\tilde z_1) \tilde \psi_1(\tilde z_1)\|_{X(H_R)}\le 2\pi
\| \tilde \psi_1\| _{X(H_R)} \sup_{\DD^1}|a_m|/(1-R'/R)$.
Finally,  we may write
$$
\EE(1/\tilde z_1, \rho z_2)\tilde  \psi_1(\tilde z_1)=
\sum_{m=0}^\infty  \rho^m z_2^m a_m(1/\tilde z_1)\tilde \psi_1(\tilde z_1) \, . 
\hbox{\qed}
$$
\enddemo

\demo{Proof of Lemma 3.2}
For $k\in \SS_0$, recall that 
$\lambda_{k,u}= \partial_1\phi_{kk,u}(0, 0) $, $ \lambda_{k,s}= \partial_2\phi_{kk,s}(0, 0)$ and
set
$$
\Delta_k= 
\cases 
\sup_{w_1 \in\partial \DD^1_k, z_2 \in \DD^2_k}
|\lambda_{k,u}  -\phi_{kk,u}(w_1, z_2)/w_1| \, , &\text{in case (P.a),}\cr
\sup_{w_1 \in -\DD^1_k, z_2  \in  \DD^2_k}
\max(|\lambda_{k,s} -\phi_{kk,s}(w_1, z_2)/z_2|,\cr
\qquad\qquad\qquad\qquad\qquad
|\lambda_{k,s}-\partial_2 \phi_{kk,s} (w_1, z_2)|)\, , &\text{in case (P.b).}\cr
\endcases
$$
We can also introduce $\Delta_k(\epsilon)$ replacing the $\DD^i_k$ by the
modified domains as when defining $\BB_0(\epsilon)$.
We shall  prove that there is $C> 0$  such that
(setting $m_k=0$ in case (P.b))
$$
\|\widehat \LL_0-\LL_0^\otimes\|_{\BB_0} \le C 
\max_{k \in \SS_0} \biggl \{
\sup (|\EE_{k, i}| +|\EE'_{k, i}|)\cdot(| \widetilde R_k-m_k|^{-1}+\Delta_k)(1- |\lambda_{k,u/s}|)^{-1}\, .
\biggr \} \tag{3.4}
$$
This will immediately imply that $\widehat \LL_0$ is bounded on $\BB_0$.

The proof of \thetag{3.4} also gives the upper bound 
$C  \max_{k\in \SS_0}\{(| \widetilde R_k-m_k|+1/\epsilon)^{-1}+\Delta_k(\epsilon)\}$
for $\|\widehat \LL_0-\LL_0^\otimes\|$ on $\BB_0(\epsilon)$.

For the spectral claim, we use the fact that
$\sup_{k \in \SS_0}\Delta_k(\epsilon)$ tends to zero as $\epsilon$ goes
to zero. We cannot apply ordinary perturbation theory,
since the Banach spaces vary. However, we can invoke Lemma 3.1 \therosteritem{2} together with
$$
(1-z \LL_0)^{-1} =
\sum_{j=0}^\infty \bigl (z (1-z\LL_0^\otimes) ^{-1}
 \circ (\LL_0-\LL_0^\otimes) \bigr)^j
\cdot (1-z\LL_0^\otimes) ^{-1} \, , \hbox{ on } \BB_0(\epsilon)\, .
$$

\smallskip
Let us prove \thetag{3.4},  considering first case (P.a). 
We concentrate on the Laplace component of the norm,
the supremum component is easier to handle. Since we may
rewrite the direct product \thetag{3.1} using
$\tilde \psi(w_1,\tilde z_2+1)$ in lieu of the Laplace transform, we have
from \thetag{2.22}
$$
\eqalign
{
s_{\phi'_{kk,s}} &(\widetilde \LL_{kk}- \widetilde \LL^\otimes_{kk}) 
(\tilde\psi) (z_1, \tilde z_2)=E_A + E_B + E_C\cr &=
\oint_{\partial D^1_k}\,  {dw_1\over 2i\pi}
\biggl (1-
{\EE_{k, 2}(w_1,1/\tilde z_2)+\tilde z_2^{-1} \EE_{k, 2}'(w_1,1/\tilde z_2)\over \tilde 
z_2 ^2} \biggr ) \cr
&\qquad\qquad\qquad\cdot
{
\tilde \psi(w_1,\tilde z_2+1+\tilde z_2^{-1} \EE_{k, 2}(w_1,1/\tilde z_2))
-\tilde \psi(w_1,\tilde z_2+1)
\over z_1 - \phi_{kk,u} (w_1, \FF^{-1}(\tilde z_2))}\cr
&-
\oint_{\partial D^1_k}\,  {dw_1\over 2i\pi}
{\EE_{k, 2}(w_1,1/\tilde z_2)+\tilde z_2^{-1} \EE_{k,2}'(w_1,1/\tilde z_2)\over \tilde 
z_2 ^2}
\biggl( {\tilde \psi(w_1,\tilde z_2+1)
\over z_1 - \phi_{kk,u} (w_1, \FF^{-1}(\tilde z_2))} \biggr )\cr
&+
\oint_{\partial D^1_k}\,  {dw_1\over 2i\pi}
{\tilde \psi(w_1,\tilde z_2+1)\over  z_1 - \phi_{kk,u} (w_1, \FF^{-1}(\tilde z_2))}
 \left ( 1-
{1-  \phi_{kk,u} (w_1, \FF^{-1}(\tilde z_2) )/z_1\over 1 - \lambda_{k,u} w_1/z_1 }\right ) \, .
}
$$
Let us bound the three terms $E_A$, $E_B$, $E_C$, taking (as we may,
by linearity and the tensor product topology)
$\tilde \psi(w_1,\tilde w_2)=\tilde\psi_1(w_1) \tilde\psi_2(\tilde w_2)$
with $\|\tilde\psi_1 \|=\|\tilde\psi_2\|=1$.

To prove $\|E_C\| \le \Delta_k/(1-|\lambda_{k,u}|)$, use on the one hand that
$|1-\lambda_{k,u} w_1/z_1| > 1-|\lambda_{k,u}|$ (since $|z_1|> |w_1|$) and
$$
 \left ( 1-
{1-  \phi_{kk,u} (w_1, \FF^{-1}(\tilde z_2)) /z_1\over 1 - \lambda_{k,u} w_1/z_1} \right )
= \left ( 
{  \phi_{kk,u} (w_1, \FF^{-1}(\tilde z_2)) /z_1 -  \lambda_{k,u} w_1/z_1
\over 1 - \lambda_{k,u} w_1/z_1} \right )
 \, .
$$ 
On the other hand
$$
\psi_1 \mapsto \oint_{\partial D^1_k}\,  {dw_1\over 2i\pi}
\biggl( { \psi_1(w_1)\over z_1 - \phi_{kk,u} (w_1, \FF^{-1}(\tilde z_2))} \biggr )
$$
is bounded on $\AA(\overline \complex \setminus \DD^1_k)$, 
uniformly as $\epsilon \to 0$, i.e., $\widetilde R_k \to \infty$ and the diameter
of $\DD^1_k$ tends to zero. To show this, note that for small $\epsilon$ 
we may perform the path integral over $\partial \DD^1_k$ and
use $\phi_{kk,u}^{-1}$ (despite the caveat in Remark 3.6 below).
Combining the above two facts with Sublemma 3.4 gives the bound on $E_C$.

Now, to estimate $E_B$, we use Sublemma 3.4 again, with
$\rho \sim |\lambda_{k,u}|$, and Sublemma~3.3 ~(1) to see that 
$$
\psi_2(\widetilde z_2 + 1) \cdot 
{\EE_{k, 2}(w_1,1/\tilde z_2)+\tilde z_2^{-1} \EE_{k,2}'(w_1,1/\tilde z_2)\over \tilde z_2 ^2}
\cdot \psi_1 (\phi_{kk,u}^{-1}(z_1))
$$
has norm bounded by $C \sup (|\EE_{k, 2}| +|\EE'_{k, 2}|)\widetilde R_k^{-2}$.

To see  that $\|E_A\|\le \widetilde R_k^{-1}$, write
$$
\tilde \psi_2(\tilde z_2+1+O_{w_1}(|\tilde z_2|^{-1})) 
-\tilde \psi_2(\tilde z_2+1)=
\tilde \psi'_2(\tilde z_2+1) \cdot O_{w_1}(|\tilde z_2|^{-1})
+\tilde \psi''_2(\Phi_{w_1}(\tilde z_2))
\cdot O_{w_1}(|\tilde z_2|^{-2}) \, .\tag{3.4}
$$
Since $\tilde z_2+1$, $\tilde z_2+1
+O_{w_1}(|\tilde z_2|^{-1})$, and thus $\Phi_{w_1}
(\tilde z_2)$, belong to a smaller
half-plane (uniformly in $w_1$ and $\tilde z_2$) we 
may apply   Sublemma 3.3 \therosteritem{2}, and, for the term involving the
second derivative, Sublemma ~3.3 \therosteritem{3}.

\smallskip

Case (P.b) is also handled by 
taking $\tilde \psi=\tilde \psi_1 \tilde \psi_2$.
For $(\tilde z_1, z_2) \in 
-H_{-\widetilde R_k}\times \DD^2_k \subset
\FF (U^1_k) \times \DD^2_k$, it is convenient to avoid 
as much as possible the simultaneous
occurrence of both variables $\tilde z_1$ and $z_2$
in the test functions, and we keep the integration over $w_2$:
$$
\eqalign
{
&s_{\phi'_{kk,s}}(\widetilde \LL_{kk}- \widetilde \LL^\otimes_{kk}) 
(\tilde\psi) (\tilde z_1,  z_2)=E'_A + E'_B + E'_C\cr
&=  
\biggl ( 1- {\EE_{k,1}(1/\tilde z_1,z_2)+
\tilde z_1^{-1} \EE_{k, 1}' (1/\tilde z_1,z_2)\over \tilde z_1 ^2}\biggr )\cr
&\, \cdot \oint_{\partial  \DD^2_k} {dw_2 \over 2i\pi}  
 {\partial_2 \phi_{kk,s} (1/\phi_{kk,u}^{-1}(1/\tilde z_1, z_2), z_2) 
\over
w_2 - \phi_{kk,s} (\phi_{kk,u}^{-1}(1/\tilde z_1, z_2), z_2)}
\cdot \bigl ( \tilde \psi_1(\phi_{kk,u}^{-1}(1/\tilde z_1, z_2))
- \tilde \psi_1(\tilde z_1-1) \bigr ) \tilde \psi_2(w_2) \cr
&-
\biggl ( {\EE_{k, 1}(1/\tilde z_1,z_2)+
\tilde z_1^{-1} \EE_{k, 1}' (1/\tilde z_1,z_2)\over \tilde z_1 ^2}\biggr )\cr
&\qquad\qquad
 \cdot \oint_{\partial   \DD^2_k}{dw_2 \over 2i\pi}  
{\partial_2 \phi_{kk,s} (1/\phi_{kk,u}^{-1}(1/\tilde z_1, z_2), z_2) 
\over w_2 - \phi_{kk,s} (\phi_{kk,u}^{-1}(1/\tilde z_1, z_2), z_2)}
\cdot \tilde \psi_1(\tilde z_1-1)\tilde \psi_2(w_2) \cr
&+
\oint_{\partial   \DD^2_k}
 {dw_2 \over 2i\pi}  
\left ( {\partial_2 \phi_{kk,s} (1/\phi_{kk,u}^{-1}(1/ \tilde z_1, z_2), z_2) 
\over
w_2 - \phi_{kk,s} ( \phi_{kk,u}^{-1}(1/ \tilde z_1, z_2), z_2)}-
{\lambda_{k,s} 
\over
w_2 - \lambda_{k,s}  z_2}
\right )
\tilde \psi_1(\tilde z_1 -1)\tilde \psi_2(w_2)\, ,
}
$$
recalling  \thetag{3.2} and \thetag{2.25--2.26}.
The norm of $E'_C$ may be bounded by a constant times
$\Delta_k$, as in case (P.a). For $E'_B$, we also adapt the argument above 
to find a bound  $C  /|\widetilde R^2_k|$. 
Finally, it is not difficult to see with
the help of Sublemma 3.3 \therosteritem{1, 4} and Sublemma~3.4
that the norm of $E'_A$ is of the order of
$ C/|\widetilde R_k|$.
\qed
\enddemo

\subhead 3.B The full symbolic transfer operator
\endsubhead

\definition{Definition 3.5 (Full symbolic operator)}
The full symbolic transfer operator 
$\widehat \LL$ is defined on the
direct sum $\BB =\bigoplus_{k \in \SS_0\cup \SS_1} 
\BB_k$  as  $\widehat \LL_0 \oplus \widehat \LL_1$ where
$$ 
\cases
\biggl ( \widehat \LL_1 (\oplus_k \psi_k) \biggr)_j = 
\oplus_{k, j \notin \SS_0\times \SS_0}
t_{kj} \cdot  \widehat \LL_{kj} \psi_k\, , j \in \SS_0 \cup \SS_1 
\,, \cr
\biggl ( \widehat \LL_0 (\oplus_k \psi_k) \biggr)_j = 
t_{jj} \cdot  \widehat \LL_{jj} \psi_j \, , j \in \SS_0 \, ,
\endcases \tag{3.5}
$$
for operators $\widehat \LL_{kj}$ 
given by \thetag{2.13, 2.14, 2.15}, respectively.
\enddefinition

\remark {Remark 3.6 (Reformulating the transfer operator)}
The jacobian $\det D \hat f_{kj} (w_1, w_2)$ 
may be expressed in terms of the pinning coordinates
as $\partial_1 \phi_{kj,u}(w_1, z_2)/  \partial_2 \phi_{kj,s} (w_1,z_2) $
where $\hat f_{kj}(w_1,w_2)=(z_1,z_2)$.
Performing two successive Cauchy residue computations (formally), we get
(note that $\partial_1 \phi_{kj,u}(w_1, z_2)$ indeed appears with a
$+$ sign, just like in \thetag{2.16}, because the $w_1$-pole is outside
of the integration curve):
$$
\bigl ( \widehat \LL (\oplus_\ell\psi_\ell) \bigr )_j (z_1, z_2)
=\sum_{k \in S} t_{kj} 
{ s_{ \phi'_{kj, s}}\over
\det D \hat f_{kj} (\hat f_{kj}^{-1}(z_1, z_2))}
\psi_k( \hat f_{kj}^{-1}(z_1, z_2))\, .\tag{3.6}
$$

The reader must beware that formula \thetag{3.6} does not
make sense in general. Indeed, in the hyperbolic case, if $(z_1, z_2) 
\in (\complex \setminus {\DD^1_j}) \times \DD^2_j$ 
then $\hat f_{kj}^{-1}(z_1, z_2)$ 
{\it may not} belong to $(\complex \setminus {\DD^1_k}) 
\times \DD^2_k$. This is connected to the fact that 
the pinning coordinates are defined for  $(w_1, z_2)
\in \DD^1_k \times \DD^2_j$ and do not necessarily extend to
$w_1 \in \complex \setminus \DD^1_k$.
However, we may always perform the $dw_2$ path integral 
as in \thetag{2.14} and in particular use  \thetag{2.16} in case (P.b) .

The sign of $\det D \hat f_{kj}$
is the product of  $s_{ \phi'_{kj, s}}$
and  the sign $s_{ \phi'_{kj, u}}$ of $\partial_1 \phi_{kj,u}$.
Note however that the real observables $\psi^\real(x,y)$ live on $I^1 \times I^2$
and that the ``real'' transfer operator $\LL^\real$ is connected to $\widehat \LL$ via
the change of coordinate
(see [Ru1, p. 1250] and [Ru4, p. 302])
$$
\psi (w_1, w_2) = \int_{I^1_k} dx {\psi^\real(x,w_2) \over w_1 - x} \, .
$$
The ``missing'' sign of $\partial_1 \phi_{kj, u}$ appears (morally) when replacing 
$z_1$ by its inverse image and we get the expected $1/|\det D \hat f_{kj}|$
factor.
\endremark

Although Remark~3.6 means that there is no simple way to express  
 \thetag{3.5} as a weighted composition operator
in general, we can in some sense pretend it is possible:
the kernel  expression for the iterates of 
$\widehat \LL$ guessed  from \thetag{3.6} by applying
the usual composition and multiplication
scheme  actually holds:

\proclaim{Lemma 3.7 (Naturality of the  transfer operator)}
For  $n\ge 1$, the difference 
$\widehat \LL^n-\widehat \LL_0^n$ acts on
$\BB$ according to
$$
\eqalign{
\bigl ( (\widehat \LL^n \
-\widehat \LL_0^n) \oplus_\ell \psi_\ell \bigr )_j 
= 
\oplus_{\vec k\in \SS^{n}\setminus \SS_0^n}
\prod_{m=1}^{n-1} t_{k_{m-1} k_m}\cdot  t_{k_{n}j}  \cdot
(\widehat \LL^n_{\vec k j} \psi_{k_1} )\, ,
}\tag{3.7}
$$
where, for any admissible $\vec kj \in \SS^{n+1} \setminus \SS_0^{n+1}$,
$$
\eqalign{
\widehat \LL^n_{\vec k j} \psi(z_1,z_2)&=
\oint_{\partial  \DD^1_{k_1}} \oint_{\partial  \DD^2_{k_1}}
{dw_1 \over 2i\pi} {dw_2 \over 2i\pi}  
 \,
 { s_{(\phi^n)'_{\vec kj}} \partial_2 \phi^{(n)}_{\vec k j,s} (w_1, z_2)  \over
w_2 - \phi^n_{\vec kj,s} (w_1, z_2)}
{\psi(w_1,w_2) \over z_1 - \phi^{(n)}_{\vec kj,u} (w_1, z_2)}\, ,
}\tag{3.8}
$$
if $k_1 \notin \SS_0$ and we replace \thetag{3.8} by
$\oint_{\partial \Gamma^1_{k_1}} \oint_{\partial \DD^2_{k_1}}$ 
or $\oint_{\partial \DD^1_{k_1}} \oint_{\partial \Gamma^2_{k_1}}$ if $k_1 \in 
\SS_0$.
Here, $s_{(\phi^n)'_{\vec kj}}$ is the sign of 
$ \partial_2(\phi^{(n)}_{\vec kj,s})$
on $\II^1_{k_1}\times\II^2_{j}$.
\endproclaim

\demo{Proof of Lemma 3.7}
To obtain the formula for the kernel, we shall
use Cauchy's theorem again. If $\vec k \in \SS_1^n$, we can follow
exactly Rugh's argument [Ru1], that we repeat for the convenience
of the reader (and because it will be adapted to $\vec k \notin
\SS_1^n$).  Denote 
$$
\GG_{\vec kj}^{(n)} (w,z) =
 { s_{(\phi^n)'_{\vec kj}} \partial_2 \phi^{(n)}_{\vec k j,s} (w_1, z_2)  \over
w_2 - \phi^{(n)}_{\vec kj,s} (w_1, z_2)}
{1
\over z_1 - \phi^{(n)}_{\vec kj,u} (w_1, z_2)}\, .
$$
By definition, $\GG_{kj}^{(1)}(w,z)$ is the kernel of 
$\widehat \LL_{ kj}$. Lemma 3.7 is clearly
true for $n=1$ and it suffices to prove inductively that for all $n\ge 1$
$$
\GG_{\vec kj}^{(n+1)} (w,z) =
\oint_{ \partial \DD^1_k} \oint_{ \partial \DD^2_k}
{d\xi_1 \over 2i\pi} {d\xi_2 \over 2i\pi}
\GG_{\vec k k_n}^{(n)} (w,\xi) \GG_{k_n j}^{(1)} (\xi,z)
\, .\tag{3.9}
$$
To prove the above equality, recall the fixed point $\xi^*=(\xi_1^*, \xi_2^*)$
constructed in the hyperbolic part of the proof of Proposition ~2.6.
The right-hand side of \thetag{3.9} has a single 
simple pole in each
coordinate at  $\xi^*=(\xi_1^*, \xi_2^*)$. Thus, writing $\text{R}_{\xi^*}$
for the residue at $\xi^*$
$$
\eqalign{
&\oint_{\partial  \DD^1_k} \oint_{\partial  \DD^2_k}
{d\xi_1 \over 2i\pi} {d\xi_2 \over 2i\pi}
\GG_{\vec k k_n}^{(n)} (w,\xi) \GG_{k_n j}^{(1)} (\xi,z)
= \cr
&\qquad \qquad \qquad \quad s_{(\phi)'_{k_n j}} s_{(\phi^n)'_{\vec k k_n}}  \text{R}_{\xi^*} \biggl [
 ( ( \xi_2 -\phi_{k_n j, s} (\xi_1,z_2))( \xi_1 -\phi^{(n)}_{\vec k k_n, u} (w_1,\xi_2))
  )^{-1} \biggr ]  \cr
&\qquad \qquad\qquad\qquad\qquad\qquad\qquad \qquad \cdot
{ 
\partial_2 \phi_{ k_n j ,s} (\xi_1^*, z_2) 
\over z_1 - \phi_{k_n j,u} (\xi_1^*, z_2)}
{ 
\partial_2 \phi^{(n)}_{\vec k k_n ,s} (w_1, \xi_2^*)  \over
w_2 - \phi^{(n)}_{\vec k k_n,s} (w_1, \xi_2^*)}
 \, .}
$$
Now, the two-variable residue is 
$ (1 -\partial_1 \phi_{k_n j,s} \partial_2 \phi_{\vec k k_n, u}^{(n)})^{-1} \, ,
$
and the fixed point property implies 
$\xi_2^*=\phi_{k_nj,s}(\phi_{\vec k k_n,u}^{(n)}(w_1,\xi_2^*), z_2)$,
so that
$$
{\partial \xi_2^* \over \partial z_2} =
\partial_1\phi_{k_n j, s}\partial_2 \phi_{\vec k k_n,u}^{(n)}
{d\xi_2^*\over dz_2} +\partial_2\phi_{k_n j, s} \, .
$$
Finally, the definition of $\phi_{\vec kj, s}^{(n+1)}$ gives
$$
\partial_2 \phi_{\vec k j, s}^{(n+1)}=
\partial_2 \phi_{ \vec k k_n, s}^{(n)} {d\xi_2^* \over dz_2}
=\partial_2 \phi^{(n)}_{\vec kk_n, s} (1-\partial_1 \phi_{k_n j, s} \partial_2 \phi_{\vec k k_n, u}^{(n)})^{-1}
\partial_2 \phi_{k_n j, s }\, .
$$
To finish the proof in the hyperbolic case,  use the multiplicative
properties of $s'_{(\phi^{n+1})_s}$.
\smallskip
If $\vec k \notin \SS_1^n$ but there are no consecutive
symbols in $\SS_0$, the above argument applies, up
to replacing $\partial \DD_k^{1,2}$ by $\Gamma_k^{1,2}$. 

So let us assume that there are at least two consecutive
$\SS_0$s.  We proceed inductively on the number
of consecutive $\SS_0$ factors, considering
the first time when $k_T , k_{T+1}\in \SS_0$ and $k_{T-1}$ or $k_{T+2} \in 
\SS_1$. 
There are four cases to consider, depending on
whether (H) is followed or preceded by (P), and on whether we are in case
(P.a) or (P.b).  In case (P.a),  we use
formula \thetag{2.14} for  
$\widehat \LL_{kk}$, and Rugh's proof for the hyperbolic case
recalled above gives the claim, using the nonhyperbolic
case of Proposition 2.6. In case (P.b), we  use \thetag{2.15}  
and the above proof may be adapted again,
using Proposition 2.6. 
\qed
\enddemo
\smallskip

The following lemma says that
if $z^{-1} \notin \sp(\widehat \LL_0)$ then the ``regularised '' transfer 
operator
$z\widehat \LL_1 (1-z\widehat \LL_0)^{-1}$ is nuclear on $\BB$ in the sense of
Grothendieck, and that its Fredholm determinant $\hat d(z)$
is dynamically defined. 
This will allow
us to describe a nontrivial part of the spectrum of $\widehat \LL=\widehat \LL_0 + \widehat \LL_1$
via this {\it regularised determinant} $\hat d(z)$,
and prove our main theorem.
Although we shall mainly refer to Grothendieck's works
[Gr1, Gr2], we mention two useful basic references: The recent book [GGK] provides a good introduction
to the theory of nuclear operator on Banach spaces, and the
survey [Ma] contains useful  results for Banach spaces of holomorphic functions.

\proclaim{Lemma 3.8
(Nuclearity of the hyperbolic analytic transfer operator)}
If $z^{-1} \notin \sp(\widehat \LL_0)$,
the operator $z\widehat \LL_1 (1-z\widehat \LL_0)^{-1}$ acting on  $\BB$
is nuclear of order zero.
The
Fredholm determinant $\det(1-z \widehat \LL_1 (1-z\widehat \LL_0)^{-1})$
is holomorphic in $\{ z^{-1} \notin \sp(\widehat \LL_0)\}$ and
$$
\hat d(z) =\det(1-z \widehat \LL_1 (1-z\widehat \LL_0)^{-1})=
\exp -\sum_{m=1}^\infty {z^m \over m}
\sum_{x \in \scriptstyle {\Fix_h} \hat f^{(m)} }  { 1 
\over  |\det (D \hat f^{(m)}_{\vec \imath_x} (x) -\Id) | }\, .\tag{3.10}
$$
\endproclaim

(Recall the bijection between hyperbolic fixed
points $x_{\vec \imath}$
of $\hat f^n$ and periodic cycles $\vec \imath_x\in \SS^{n}\setminus
\SS_0^{n}$ from Corollary~2.8.)

Since $\hat f$ is invertible the matrix equality 
$
\det (1+A)\det (1+B) = \det((1+A)(1+B))
$ gives
$$
\hat d(z)
= 
\exp -\sum_{m=1}^\infty {z^m \over m}
\sum_{x \in \scriptstyle {\Fix_h} \hat f^{(m)}} { 1\over
 |\det D \hat f^{(m)}_{\vec \imath_x} (x)|} {1\over
  |\det ( \Id- D \hat f^{(-m)}_{\vec \imath_x} (x)) | }\, .\tag{3.11}
$$
\medskip

\demo {\bf Proof of Lemma~3.8}
We adapt the argument in [Ru3, Lemma 2.7].

Using [Gr1], one shows that for any compact sets $K$, $K'$ of the
complex plane such that $K$ is contained in the
interior of $K'$, the restriction map from the
Banach space $\AA(K')$ of analytic functions on $\Int K'$ extending continuously
to $K'$, to the space $\AA(K)$ of analytic 
functions on $\Int K$ extending continuously to $K$,  is nuclear of order zero.  (See e.g. 
the beginning
of the proof of Lemma~2.7 in [Ru3], which was adapted from Ruelle's paper
[Rue1].) This argument may be extended to show that if
$K^1$, $K^2$, $(K^1)'$ and $(K^2)'$ are compact subsets of
$\complex$ such that
$$
K^2\subset \Int (K^2)' \, , \qquad (K^1)'\subset \Int  K^1 \, ,\tag{3.12}
$$
 then the restriction map $r_{K}$ from the Banach space 
$\AA(\overline \complex \setminus (K^1)', (K^2)')$ 
of functions holomorphic in the interior of
$(\overline \complex \setminus (K^1)')\times (K^2)'$, vanishing at infinity,
and extending continuously to the boundary,
to the Banach space $\AA(\overline \complex\setminus K^1, K^2)$,  is nuclear
of order zero. (The key step is to observe that
the topological vector space of analytic functions on  
$((\overline {\complex} \setminus K^1)\times \Int\,  K^2)$
is a nuclear space.)

On $\BB_1 =\oplus_{k\in \SS_1}   \BB_k$, the resolvent
$(1-z\widehat \LL_0)^{-1}$ acts as the identity,
while $(1-z\widehat \LL_0)^{-1}$  is bounded on $\BB_0$
by our assumption on $z$.  Thus, $(1-z\widehat \LL_0)^{-1}$ is bounded on $\BB$.

Recall the definitions of $U^2_k$ and $U^1_k$ for $k\in \SS_0$  from Section~2.
Let $(K^1_k)', (K^2_k)'$ be compact sets
such that $(K^2_k)'\subset \Int \DD^2_k$, if
$k\in \SS_1$ or we are in case (P.b), and
$(K^2_k)'\subset U^2_k$ in case (P.a), while   $\DD^1_k\subset \Int \, 
(K^1_k)'$
if $k\in \SS_1$ or we are in case (P.a),
and   $(K^1_k)'\subset U^1_k$ in case (P.b). 
Denote 
$$
\BB(K'):=\biggl [ \oplus_{k\notin (P.b))}
 \AA(\bar \complex \setminus (K^1_k)',(K^2_k)')\biggr ]
\oplus
\biggl [\oplus_{k \in (P.b)} \AA((K^1_k)',(K^2_k)') \biggr ]\,.
$$
The restriction map $r_{K'} : \BB \to  \BB(K')$
is continuous (see [Ru3, Lemma 2.6] for case (P.a)).
Therefore $r_{K'}\circ  (1-z\widehat \LL_0)^{-1}$ is 
continuous from $\BB$ to $\BB(K')$.

Recall the definition of the compact set  $K^2_k$ for $k \in \SS_0$. Extend 
it 
to $k \in \SS_1$ as follows:
$$
K^2_k=\cup_{\ell  \in \SS} \phi_{k\ell ,s} (\DD^1_k,  \DD^2_\ell)\, .
$$
For $k$ in  case (P.b) we take a compact subset $K^1_k$ 
of $ \Int (K^1_k)'\subset U^1_k$. In cases (H) or (P.a) we take a compact 
subset $K^1_k$ of
$\widetilde \DD^1_k$, containing $\DD^1_k$ in its interior and so that
$\phi_{kj, u}(K^1_k, \DD^2_j) \subset \Int \DD^1_j$ for each  $j$.
Now, up to slightly changing the $(K^{1,2}_k)'$ introduced above, we may 
ensure that
$K^2_k\subset \Int (K^2_k)' $ and $ (K^1_k)'\subset \Int  K^1_k$ 
in cases (H) or (P.a), 
$K^2_k \subset \Int (K^2_k)'$ and $K^1_k \subset \Int (K^1_k)'$ in case (P.b), while 
maintaining the other requirements. By the above choices, the restriction $r_K$
is nuclear from $\BB(K')$ to $\BB(K)$.

Then (use in particular the definition of $\Gamma^{1;2}_k$ if
$k \in \SS_0$)
$\widehat \LL_1$ is bounded from 
$$
\BB(K)=\biggl [\oplus_{k\notin (P.b)} \AA(\bar \complex \setminus K^1_k,K_k^2)\biggr ]
\oplus
\biggl [ \oplus_{k \in (P.b)} \AA(K^1_k,K^2_k)\biggr ]
$$
to  $\BB':=\oplus_{k \in \SS}\BB'_k$.
Since the inclusion $j : \BB' \subset \BB$
is continuous (for $k\in \SS_0$ see [Ru3, Lemma 2.5] for a similar result, noting that
the constant $\beta$ there is equal to $+1$ in our case), 
the composition  $j \circ \widehat \LL_1$,  
is bounded from $\BB(K)$ to $\BB$. 

As a consequence of the above considerations, the composition
$$
\widehat \LL_1  (1-z\widehat \LL_0)^{-1} 
=\biggl ( j\circ  \widehat \LL_1\biggr )  \circ r_{K} 
\biggl (\circ r_{K'} \circ  (1-z\widehat \LL_0)^{-1} \biggr ) : \BB\to \BB
$$
is nuclear  of order zero (just use that a nuclear
operator composed with  bounded operators is nuclear). 
By [Gr1, II, pp.~16 and~18]
it has a Fredholm determinant 
$\det(1- z\widehat \LL_1  (1-z\widehat \LL_0)^{-1})$ which is an entire
function of $1/z \notin \sp(\widehat \LL_0)$.

\smallskip
It remains to establish the stated ``dynamical'' formula for the 
traces. For this, we combine arguments
from [Ru3] and [Ru2].
First, just like on p.~ 17 of [Ru3], we find that
for small enough $z$ 
$$
\tr \log (1-z \widehat \LL_1 (1-z\widehat \LL_0)^{-1})
= \sum_{m=1}^\infty
{z^m \over m} \tr (( \widehat \LL_0 + \widehat \LL_1)^m -\widehat \LL_0^m) \, .
$$ 
Defining $d_m =\tr (\widehat \LL^m -\widehat \LL_0^m)$,
the Fredholm determinant
$$
\det (1-z \widehat \LL_1 (1-z\widehat \LL_0)^{-1})=
\exp- \sum_{m=1}^\infty {z^m\over m} d_m 
$$
extends analytically to $\complex - \sp(\widehat \LL_0)$.
Using the uniform contraction, we get for the trace
$$
d_m = \sum _{\vec \jmath \in \SS^{m+1} \setminus \SS_0^{m+1}, j_1=j_{m+1}}
\prod_{k=1}^{m} t_{j_{k} j_{k+1}}  \cdot
\tr (\widehat \LL^m_{\vec \jmath}  )\, ,
$$
where the iterated operator $\widehat \LL^m_{\vec \jmath}$
may be expressed in kernel form by \thetag{3.8}. 
Finally, the trace may be computed by performing
a Cauchy integration:
$$
\eqalign{
\tr (\widehat \LL^m_{\vec \jmath })&=
\oint_{\partial \DD^1_{j_1}} \oint_{\partial \DD^2_{j_1}}
{dw_1 \over 2i\pi} {dw_2 \over 2i\pi}  
 { s_{(\phi^m)'_{\vec \jmath}} 
\partial_2 \phi^{(m)}_{\vec \jmath,s} (w_1, w_2)  \over
w_2 - \phi^m_{\vec \jmath,s} (w_1, w_2)}
{1 \over w_1 - \phi^{(m)}_{\vec \jmath,u} (w_1, w_2)}\, \cr
&
= {1 \over |\det (D\hat f^m(x_{\vec \jmath}) -\Id)|}\, ,\cr
}
$$
where $x_{\vec \jmath}$ is the unique (necessarily
hyperbolic, and real) fixed point of 
$\hat f^m$ associated to the admissible periodic sequence
$\vec \jmath \in \SS^{m+1} \setminus \SS_0^{m+1}$ (see pp. 1246--1247 of [Ru1]
for details). \qed
\enddemo

To relate the zeroes of the (regularised)
Fredholm determinant $\hat d(z)$ 
(to part of) the spectrum of $\widehat \LL$ we need
the following lemma of Rugh, that we state for convenience of the reader:

\proclaim{Lemma [Ru3, Lemma 2.8]}
Let $\MM_0 : \BB \to \BB$ be a bounded linear operator on a Banach space
$\BB$ and let $\MM_1 : \BB \to \BB$ be nuclear of order zero.
Assume that $\Bbb C\setminus \sp(\MM_0)$ is connected. Then 
the part of the spectrum of $\MM_0+\MM_1$ which does
not intersect $\sp(\MM_0)$ consists of isolated 
eigenvalues of finite multiplicity,
which can not accumulate in  $\Bbb C \setminus\sp(\MM_0)$.
The Fredholm determinant
$$
d(u) =\det(1-\MM_1(u-\MM_0)^{-1})
$$
is analytic in $u \in \Bbb\complex \setminus \sp(\MM_0)$. In this domain,
the zero-set of $d(u)$ counted with order is the same as
the eigenvalues of $\MM_0+\MM_1$ counted with (algebraic)
multiplicity.
\endproclaim

We do not know a priori that the complement of the
spectrum of $\widehat \LL_0$ is connected. However,
the spectra of the direct product operators \thetag{3.1} and \thetag{3.2} 
have this property  by Lemma 3.1.  We saw in Lemma~3.2
how to compare $\widehat \LL_0$  to $\widehat \LL_0^\otimes$,
which is a direct sum of such operators.
In our application (see Section~ 4) we will find a closed set containing the
spectrum of $\widehat \LL_0$, arbitrarily close to
the spectrum $\widehat \LL_0^\otimes$, with connected complement
in $\Bbb C$, and apply the above lemma to this complement.

%%%%%%%%%%%%%%%%%%%%%%%%%%%%%%%%%%%%%%%%%%%%%%%%%%%%%%%%%%%%%%%%%%%%%

\head 4.  Reducing to (symbolic) analytic almost hyperbolic maps
\endhead

Let us consider now a real-analytic diffeomorphism $f: M \to M$
for which there exists a dominated splitting $T_\Omega = E\oplus
F$ over the nonwandering set. Our starting point
will be a  decomposition of $\Omega$ from [PS2].
The decomposition in [PS2] is stated for the limit
set of $f$. However, if $M$ is a surface and $\Omega$ 
is hyperbolic then $\Omega$ coincides with the limit set of $f$ (see [NP]), and 
this equality also holds  when $\Omega$ has a
dominated splitting.

We recall some notation and results from [PS2]. We say that a
compact invariant set $\Lambda \subset \Omega$ admits a {\it
spectral decomposition} if it is a finite disjoint union of
transitive compact invariant sets $\Lambda_i$ (called {\it basic sets}) which
may further be decomposed in a finite union of $n_i\ge 1$ {\it basic subsets}
$\Lambda_{i_j}$ with $f(\Lambda_{i_j}) =\Lambda_{i_{(j+1 \mod
n_i)}}$, and $f^{n_i} |_{\Lambda_{i_j}}$  topologically mixing. 
We shall assume that the $\Lambda_j$ are not trivial, i.e.,
not reduced to a single periodic orbit. It
follows from the results in [PS1] and the classical
Hirsch-Pugh-Shub [HPS] theory that for each small enough
$\epsilon$, there is $\delta$, so that for each $x \in \Omega$, there
exist local center stable and unstable manifolds
$W^{cs}_\epsilon(x)$ and $W^{cu}_\epsilon(x)$ so that
$$
\eqalign { &T_x W_\epsilon^{cs}(x)=E \, , \quad T_x W_\epsilon^{cu}(x)=F
\, ,\cr &f(W_\delta^{cs}(x)) \subset W^{cs}_\epsilon(f(x))\, , \quad
f^{-1}(W_\delta^{cu}(x)) \subset W^{cu}_\epsilon(f^{-1}(x))\, . }
$$

The decomposition proved by Pujals and Sambarino [PS1] (who do not require analyticity,
$C^2$ suffices)
says that $\Omega = \Lambda \cup \RR\cup \II$.
Here, the ``quasi-periodic'' set
$\RR$ is a finite
union of normally hyperbolic $C^2$ simple closed curves $\CC_i$ on which $f^{r_i}$ is conjugated 
to an irrational rotation for $r_i\ge 1$. 
The ``periodic'' set $\II$ is  the
union of a finite set of isolated periodic orbits with a set
contained in  a finite union $\cup_j \II_j$ of normally
hyperbolic  $C^2$ arcs or simple closed curves with $f^{m_j}(\II_j)\subset \II_j$
for $m_j\ge 1$. The set $\II$  contains all  $\Omega\setminus P$
isolated periodic orbits.
Next,
$f$ is expansive on the ``almost hyperbolic'' compact invariant
set $\Lambda$, which admits a  spectral
decomposition  $\Lambda = \cup_j \Lambda_j$, together
with local product structure. (Lemma 4.5.1 in [PS2]: there are $\gamma$ and $\eta > 0$ so that for
any $x,y\in \Lambda_j$ with $d(x,y) < \eta$ then $W^{cs}_\gamma(x)
\cap W^{cu}_\gamma(y) \in \Lambda_j$.) 
Finally,  the set $\NN$ of nonhyperbolic periodic orbits in $\Lambda$ 
is empty or finite.
In fact (see e.g. Proposition~A.2 below), every basic set $\Lambda_j$
which does not contain any nonhyperbolic periodic point is
uniformly hyperbolic.

It is easy to construct examples where $\cup \II_j$ is not empty:
just take a real-analytic flow on the sphere with
both poles as sources and the equator as limit set.
Our analyticity assumption implies that the arcs and curves $\II_j$ in
$\II$ are isolated: Indeed, if they were not isolated, there would be a
normally hyperbolic arc $\II_j$ and a basic set $\Lambda_k$
such that their intersection is a nonhyperbolic periodic point
$q$ which is accumulated by periodic points contained in
$\II_i$ with the same period than that of $q$; and this would
contradict analyticity. It follows that $\II$ is not
only open but also compact in $\Omega$. 
Note also that the set $\HH$ of  (isolated) hyperbolic  periodic points in $\II$ is finite.
(Indeed, if there were infinitely many hyperbolic periodic points
in $\cup_j \II_j$, their periods being bounded by $\max m_j$, a subset of
constant period  would accumulate on a periodic point, contradicting
the analyticity assumption.)

Note that $\RR$ does not contain any periodic orbits.
Consider first the finite set $\HH \subset \II$ of isolated hyperbolic periodic orbits. 
Writing  $P\ge 1$ for the
period and $\lambda_E$, $\lambda_F$ for the multipliers (eigenvalues of
$Df^P(p)$), of a periodic orbit $p$, each
$p\in \HH$ contributes to $d_f(z)$ a factor of the following type:
$$
\eqalign{
d_{f|{\text sink}} (z)&= \prod_{j=0}^{\infty} \prod_{k=0}^\infty
(1-z^P\lambda_E^{j}\lambda_F^{k})\, , \quad |\lambda_{E}|<  |\lambda_F|< 1 \, , \cr
d_{f|{\text saddle}} (z)&= \prod_{j=0}^{\infty} \prod_{k=1}^\infty
(1-z^P\lambda_E^{j}\lambda_F^{-k})\, , \quad |\lambda_{E}| < 1\, , |\lambda_F|> 1 \, , \cr
d_{f|{\text source}} (z)&= \prod_{j=1}^{\infty} \prod_{k=1}^\infty
(1-z^P\lambda_E^{-j}\lambda_F^{-k})\, , \quad 1< |\lambda_{E}| < |\lambda_F| \, . \cr
}
\tag{4.1}
$$
The infinite products above  all converge, and define
entire functions with an obvious zero-set.
In particular, each $d_{f|{\text sink}} (z)$ is zero-free in the
open unit disk and admits $P$ simple zeros on the closed disk,
at the $P$th roots of $1$, while each $d_{f|{\text source}} (z)$
admits a first zero at $z^P=\lambda_E\lambda_F$, which is
outside the open disk, and each $d_{f|{\text saddle}} (z)$
admits a first zero at $z^P=\lambda_F$, which is
outside the open disk.

We may therefore concentrate on the dynamical determinant
$\prod_j d_{f|\Lambda_j}(z)$, where the $\Lambda_j$ are the basic
sets of $\Lambda$. 
Recall the set $\Sigma(p)$ associated to $p\in \NN$ by \thetag{1.3}.
To prove Theorem~A, we need to see that
$d_{f|_{\Lambda_j}}(z)$ is holomorphic in the (possibly)
slit plane, or multiply slit plane defined by
$
\{ z \in \complex \mid  1/z \notin \cup_{p\in \NN\cap \Lambda_j}\Sigma(p) \}
$.

In order to do this, we shall associate an almost hyperbolic
analytic map $\hat f$ to $f|{\Lambda_j}$ in such a way as
to ensure that $d_{\hat f}(z)$ is almost $d_{f|{\Lambda_j}}(z)$
(dealing with the usual overcounting of periodic orbits on
the boundaries of the Markov rectangles is postponed until Section~5.A).

\smallskip

{\bf Markov partitions.}

The starting point in our construction of the symbolic map $\hat f$ is the 
existence of Markov partitions for $f$. We recall a possible definition in dimension two:

\definition{Definition 4.1 (Markov partition)}
Let $\Lambda \subset  \Omega$ be a 
basic set. A Markov partition $\RR$ of $\Lambda$ is a finite
collection $\{R_1, \ldots, R_\ell\}$ of ``rectangles,'' with disjoint
interiors, which are diffeomorphic to the square $Q=[-1,1]^2$, through
$R_i=\tilde \psi_i(Q)$, whose union contains $\Lambda$,  and such that
$$
f(\partial_s R_i) \subset \cup_j \partial_s R_j \, ,
f^{-1}(\partial_u R_i) \subset \cup_j \partial_u R_j \, ,
$$
where $\partial_s R_i =\partial_E R_i =\tilde \psi_i(\{ (x,y) \mid |y|=1\})$ and
$\partial_u R_i =\partial_F R_i =\tilde \psi_i(\{ (x,y) \mid |x|=1\})$.
\enddefinition

To a Markov partition with $\ell$ rectangles we may associate
an $\ell \times \ell$ transition matrix by setting
$t_{ij}=1$ if the interior of $f(R_i)$ intersects
$R_j$ and $t_{ij}=0$  otherwise. Transitivity of $\Lambda$
implies that this matrix is irreducible with no wandering states.

Since we are in dimension two, we
can adapt the construction of Markov partitions in  
 [PT, Appendix~2] (see [PS2, Lemma 4.5.2], first reduce to a mixing basic subset). 
The construction shows that  a basic set $\Lambda$ of 
$\Omega$ admits Markov partitions of arbitrarily small diameter
(the diameter being the maximum of the diameters of the rectangles
$R_i$). Since we have only a
finite number of nonhyperbolic periodic points in $\Lambda$, we may
assume that  each rectangle
contains at most one nonhyperbolic periodic point. 
We may furthermore ensure that if $q\in R_i$ is
$E$-nonhyperbolic (i.e., $\lambda_E=\pm1$) then $q\in \partial_u R_i=\partial_F R_i$ 
but $q \notin \partial_s R_i$,
while if $q$ is
$F$-nonhyperbolic then $q\in \partial_s R_i$ but $q \notin \partial_u R_i$. 

Note that if $q$ is fixed,
nonhyperbolic, and the order of $f-\Id$ at $q$ is even (i.e., we have
a saddle node situation), then $\Omega$
lies entirely on the weakly attracting side of $q$
if $\lambda_E=1$, while it is on the weakly
expanding side of $q$ if $\lambda_F=1$. If the nonhyperbolic
multiplier is $+1$ but the order of $f-\Id$ at
a nonhyperbolic fixed point is odd (i.e., we have a saddle), then $\Omega$ intersects both sides,
so that $q$ will belong to the boundary of two rectangles.
If the nonhyperbolic multiplier is $-1$ then $\Omega$ also meets both sides
and we need two rectangles,
whether the order is even or odd. 
For periodic points of period larger than one, the above remarks
may be applied along the orbit.

We set $i\in \SS_0$
if $R_i\cap \NN\ne \emptyset$ (i.e., it is a ``bad'' rectangle)
and $i\in \SS_1$ otherwise.

\smallskip

{\bf Takens ($C^\infty$) local coordinates for nonhyperbolic fixed points.}

In the arguments below it will be convenient to use normal
forms. We discuss first  the $C^\infty$ normal form due to Takens

Let $q\in \NN\cap \Lambda_j$ be  an  $F$-nonhyperbolic 
fixed point. In particular,
$q$ is not contained in a periodic curve.
We assume also that $\lambda_F=+1$.
(The other cases in $ \NN\cap \Lambda_j$, $E$-nonhyperbolic, period $\ge 2$,
multiplier $-1$,  are similar, see also \S5.B.)
By [T], we may express the diffeomorphism $f$ in $C^\infty$  local coordinates at $0\mapsto q$ as
$$
f_T(s,t)=(\lambda(t) s \,  , b(t)) \, , \tag{4.2}
$$
for $C^\infty$ functions $\lambda$ and $b$ satisfying
$0<|\lambda(0)|<1$, $b(0)=0$, and $b'(0)=1$.
(Notice that $\{(s,t): t=0\}$
is the strong stable ($E$) manifold  and $\{(s,t):s=0\}$ is the central unstable ($F$)
manifold.)

\smallskip

The neutral Takens coordinate $b(t)$ cannot be infinitely flat in our setting:

\proclaim{Lemma 4.2 (Nonflatness of $C^\infty$ normal form)}
If $f$ is analytic and $q$ is not contained in a curve
of fixed points,  setting $\nu+1\ge 2$ to be
the multiplicity of $f-\Id$ at $q$, then $b^{(\nu+1)}(0)\neq 0$,
and $b^{(j)}(0)=0$ for $2\le j \le \nu$.
\endproclaim

\demo{Proof of Lemma~4.2}
Since $q$ is not contained in a curve of fixed points,
by [H, Proposition 2.3, see also p. 481], we may express the diffeomorphism $f$ in
real-analytic local coordinates at $0$ as
$$
f_H(x,y)=(g(x) +yh(x,y)\, , y+ y^{\nu+1}+A y^{2\nu+1}+a(x)y^{2\nu+2}+\cdots) \, , \tag{4.3}
$$
for  $\nu+1\ge 2$ the multiplicity of $f-\Id$ at $0$,  with $A$ a   complex constant and
$g$, $a$, and   $h$ real-analytic in a neighbourhood of
$0$ in $\complex$ (resp. $\complex^2$) and $g(0)=0$, $h(0,0)=0$ and
$0<|g'(0)|<1$. (Hakim deals with holomorphic situation, but real-analytic
data gives real $A$ and real-analytic functions $g$, $h$, and $a$.)
Notice that in these coordinates $\{(x,y): y=0\}$
is the still strong stable manifold. The central
manifold, however, does not have an obvious description any more
(indeed, it is usually not real-analytic), but it is tangent to
$\{(x,y):x=0\}$ and can be described as the
graph of a $C^\infty$ map $y\mapsto x_F(y)$ with $x_F(0)=0$.  Additionally, it is the image of
$\{(s,t):s=0\}$ by the conjugacy restricted to this line, which
may be encoded in the $C^\infty$ one-dimensional
diffeomorphism $t\mapsto y_t$ with inverse $y\mapsto \tau(y)$ and $y_0=0$, $y'_0\ne 0$.
Now, $b(t)$ can be decomposed as
$$
t\mapsto f_H(x_F(y_t), y_t)\mapsto 
b(t)= \tau( y_t+y_t^{\nu+1}+Ay_t^{2\nu+1}+a(x_F(y_t))y_t^{2\nu+2}+\cdots) \, .
$$
Hence, using the mean value theorem, and setting
$$e(t)=y_t^{\nu+1}+Ay_t^{2\nu+1}+a(x_F(y_t))y_t^{2\nu+2}+\cdots\, , $$ 
we have
$
b(t)= t+ \tau'(\xi(t)) e(t)
$.
Since $e^{(\ell)}(0)=0$ for $0\le \ell < \nu+1$, it follows that  
$b^{(\ell)}(0=0$ for $\ell < \nu+1$ and 
$
b^{(\nu+1)}(0)=\tau'(0) e^{(\nu+1) } (0) =
\tau' (0) (\nu+1)! \, (y'_0)^{\nu+1}\neq 0
$. \qed
\enddemo

\smallskip

{\bf Constructing a symbolic model}

From now on, we work with a real-analytic
atlas of $M$, $\psi_k : A_k\to M$,
where each $A_k\subset \real^2$ is a viewed as a real subset of the complexification
$T^\complex_{x} M\subset \complex^2$ 
of the tangent space for some chosen
$x \in \hbox{Int} \psi_k(A_k)$ (we refer to 
pp.~808--809  in [Ru2] for details). 
We assume that the atlas is compatible with the normal
form $f_H$ from \thetag{4.3}. We let $\widehat A_k$ be a complex neighbourhood
of $A_k$. Note that
the decomposition $E\oplus F$ extends to  $T^\complex_{\Lambda} M=E^\complex_{\Lambda}+F^\complex_{\Lambda}$.
For $x \in A_k$ with
$\psi_k(x) \in \Lambda$ we denote $E^s_{k,x}=E^\complex_{k,x}$, $E^u_{k,x}=F^\complex_{x}$,
let $P^{u, \complex}_{k,z} :T^\complex_{z} M \to
E^{u,\complex}_{k, z}$, $P^{s, \complex}_{k,z}: T^\complex_{z} M \to E^{s,\complex}_{k, z}$
denote the complexified projections to the unstable and stable bundles
(we sometimes drop the $k$ index). 
By construction, for  $j,k$
with $t_{j,k}\ne 0$, the map $\hat f$ induced by $f$ 
in the charts extends to a real-analytic map.

We are now almost ready to state Proposition~4.5, which says that
a (sequence of) almost hyperbolic symbolic models for $f$ can be constructed with
the help of a sequence of
Markov partitions of diameters ending to zero.
(All Markov partitions involved will be real and compatible with the real-analytic 
atlas chosen above, in the sense
that each rectangle is included in some $\psi_j(A_j)$.) We must introduce
further notation:

\definition{Definition 4.3 (Admissible complex extension)} 
A subset $\widetilde \Lambda$ of
a complex neighbourhood of $\Lambda$ (in the charts) is an admissible complex extension
of $\Lambda$ if there is a complex neighbourhood $V_j$ of each nonhyperbolic fixed point
$q_j\in \NN$ such that (in the charts) $\widetilde \Lambda\setminus \cup_j V_j$
is a complex neighbourhood of $\Lambda\setminus \cup_j V^\real_j$,
and each $\widetilde \Lambda\cap V_j$ in  analytic charts $(z_E, z_F)$ 
(compatible with \thetag{4.3}) contains the intersection
of a neighbourhood of $q_j=0$ and a domain
$
\{\Re (z_E ^{\nu_j})\ge  0 \}
$,
if $q_j$ is $E$-nonhyperbolic,  and
$
\{\Re (z_F^{\nu_j})\ge  0  \}
$ otherwise (as usual, $\nu_j+1\ge 2$ is the multiplicity).
\enddefinition

Let $| \cdot |_z$ denote the norm on
$T^\complex_{z} M $ induced by Riemann metric. 
Adapting the Mather trick and ideas from Crovisier, we show in the Appendix:

\proclaim{Lemma 4.3 (Adapted metrics)}
Assume all $q\in \NN$ are fixed points with nonhyperbolic multiplier $+1$.
There are two 
semi-norms $\|\cdot \|_{E,z}=\|\cdot \|^s_{z}$, $\| \cdot\|_{F,z} =\| \cdot\|^u_{z}$
on the complex tangent bundle $T_{ \Lambda^\complex}M$ 
over a complex neighbourhood $\Lambda^\complex$ of $\Lambda$,
and 
an admissible complex extension $\widetilde \Lambda$ 
of $\Lambda$, such that:
\roster
\item
For all  $w$
$$
 \|Df^{-1} w \|_{F,f^{-1}z}\le \|w\|_{F,z}\, , \qquad \|Df w\|_{E,f z}\le \|w\|_{E,z} 
\, , \forall z \in \widetilde \Lambda \, . \tag{4.4}
$$
For  $w\ne 0$, equality holds in the first bound if  and only if
$z$ is  $F$-nonhyperbolic,
in the second one if and only if $z$ is 
$E$-nonhyperbolic. In addition,
there are $C>0$ and  $C_\nu > \nu$ so that for each
$F$-nonhyperbolic fixed point $q_j\in \Lambda$  of index $\nu+1$, 
letting $V_j$ be the neighbourhood from Definition 4.3:
$$
\|Df^{-1}w\|_{F, z}\le (1-C|z_E|)(1-C_\nu |\Re (z_F^\nu )|)  \|w\|_{F,f(z)}\, ,
\forall z=(z_E, z_F) \in \widetilde \Lambda \cap V_j\, ,  \tag{4.5}
$$ 
and similarly for the $E$-nonhyperbolic case.
\item
There is $C< \infty$ so that   
$$
{ |\cdot|_z\over C} \le 
\max(\|\cdot \|_{E,z}  ,  \| \cdot\|_{F,z}) \le  C |\cdot|_z\, , \quad
 \forall z\in  \Lambda^\complex\, .\tag{4.6}
$$
\item
For any $z$ in $\Lambda^\complex$, and any $z_0$ close enough to
$z$ so that $(z-z_0)$ can be interpreted as a vector and
$(f (z) -f(z_0) - Df_{z_0} (z-z_0))$ is in  the chart at $f(z_0)$, we have
$$
\| f (z) -f(z_0) - Df_{z_0} (z-z_0) \|_{E/F, f(z_0)}
\le \epsilon^{E/F}_1(|z-z_0|_{z_0}) |z-z_0|_{ z_0} \, ,
$$
where  the two functions $\epsilon^{E/F}_1:\real^+\to \real^+$ are  
$\CC^\infty$  with $\epsilon^{E/F}_1(0)=0$.
Furthermore, there is $C$ so that
for any $\delta$, if $z$ and $z_0$ are in a complex
$\delta$-neighbourhood of an $E/F$-nonhyperbolic fixed point $q\in \Lambda$ of index $\nu+1$, then
$$
\epsilon_1^{E/F}(|z-z_0|_{z_0}) \le C  \cdot \delta^{\nu} |z-z_0|_{z_0}\, .  \tag{4.7}
$$ 
\item Let $\gamma >0$ be the Hoelder smoothness of the stable and unstable
foliations. Then:
$$
 |\|\cdot\|_{E/F,z} - \|\cdot\|_{E/F,z'}|\le \epsilon_{2}( d(z,z')) | \cdot|_z
\, , \quad \forall z , z'\in \Lambda^\complex \, , \tag{4.8}
$$
for a $\gamma$-Hoelder continuous function $\epsilon_2$ on $\real^+$, vanishing at $0$.
If $z$ and $z'$ are on the same $W^{E/F}$ local stable manifold,
we can assume that $\gamma=1$.
\endroster
\endproclaim

We set
$
\| \cdot\|_{z} =\max(\| \cdot\|^u_{z},\| \cdot\|^s_{z}) 
$
We shall work with two types of complex extensions of the
rectangles of a Markov partition of $\Lambda$ (compatible 
with charts). If $R_k \cap \NN=\emptyset$,
for $\xi_k \in R_k \subset \psi_j(A_j)$, and small $\delta_k> 0$,
we consider an {\it $\omega$-rectangle} (just like in [Ru2])
$$
R^\omega_k (\xi_k, \delta_k) = \{ z \in \widehat A_j \subset \complex^2 \mid 
\| z-\psi^{-1}_j (\xi_k)\|_{\psi^{-1}_j (\xi_k)} \le \delta_k\} \, .
$$
By definition, $R^\omega_k$ factorises as $\DD^1_k \times \DD^2_k$, with $\DD^i_k$ 
a compact connected subset
of $\complex$ with smooth boundary (in fact, a disc), and intersecting the real axis on an interval
$\II^i_k$.

If $R_k \cap \NN=\{q\}$, we shall assume that $q=\psi_j(0)$ in charts
$z=(z_E, z_F)$ compatible with the Hakim normal form \thetag{4.3}, and,
if $q$ is $F$-nonhyperbolic of index $\nu_k+1$,  for $\delta_k>0$ 
and $\pi/(2\nu_k)< \theta_k< \pi/\nu_k $ we consider an {\it $\omega$-petal} 
$$
R^\omega_k (\nu_k, \theta_k,\delta_k)= 
\{ z \in \widehat A_j \mid z_F \in \UU(\theta_k, \delta_k)  \hbox{ and }
\| z\|^E_{q} \le \delta_k\} \, ,
$$
Again, $R^\omega_k$ factorises as $\DD^1_k \times \DD^2_k$, with $\DD^i_k$ 
a compact connected subset
of $\complex$ with smooth boundary (except at $0\in \DD^2_k$), and intersecting the real axis 
on an interval $\II^i_k$.  
If $q$ is $E$-nonhyperbolic, we proceed in an analogous
way. We denote by  $\widetilde  R^\omega_k$ the real projection
$
\widetilde  R^\omega_k = \psi_j (R^\omega_k \cap \real^2) 
$ of an $\omega$-rectangle or an $\omega$-petal.

\smallskip

We may finally state the main result of this section:

\proclaim{Proposition 4.5
(From dominated splitting to almost hyperbolic)}
Let $f$ be a real
analytic diffeomorphism on a compact real-analytic surface $M$,
with dominated splitting on its nonwandering set $\Omega$.
Let $\Lambda$ be a basic set of $f$.  Assume that all orbits in $\NN$
are fixed points with neutral multiplier $+1$ and multiplicity
$\nu+1=2$. Then there exists a sequence of 
Markov partitions $\RR_n=\{R_{k,n}\}_{k\in \SS_n}$ 
of $\Lambda$, with diameters tending to zero and such that, for each fixed $n$,
denoting by $t_{ij}=t_{ij,n}$ the transition matrix:
\roster
\item  For each   $k$ so that  $R_k \cap  \NN=\emptyset$,
letting  $j$ be such that
$R_k \subset \psi_j(A_j)$, there are  $\xi_k \in R_k$ and $\delta_k> 0$ so that 
the  projection of the corresponding $\omega$-rectangle satisfies
$R_k\subset \widetilde R^\omega_k(\xi_k, \delta_k) \subset \psi_j(A_j)$.

\item   For each $k$ so that  $R_k \cap  \NN=\{q\}$, of index $\nu_k+1=2$,
letting  $j$ be such that $R_k \subset \psi_j(A_j)$,  there are   
$\delta_k> 0$ and $\pi/(2\nu_k) < \theta_k< \pi/\nu_k$,  so that the
projection of the corresponding $\omega$-petal
satisfies
$R_k\subset \widetilde R^\omega_k(\nu_k, \theta_k, \delta_k) \subset \psi_j(A_j)$.

\item  The following defines an almost hyperbolic analytic map $\hat f$:
$$
\eqalign{
&\SS_0=\{k \mid R_k \cap \NN \ne \emptyset\} \, , \quad \SS_1=\{k \mid R_k \cap \NN =\emptyset\}
\, , \quad 
R^\omega_k=\DD^1_k \times \DD^2_k \subset \widehat A_k\, , \cr
&\hbox{for } t_{k\ell}\ne 0 :\quad\psi_i \circ \hat f_{k\ell}|_{(\Cal I^1_k \times \Cal I^2_k)
\cap \hat f_{k\ell}^{-1} (\Cal I^1_\ell\times \Cal I^2_\ell)}
= f \circ \psi_j |_{\psi_j^{-1} (\widetilde R^\omega_k \cap f^{-1} (\tilde  R^\omega_\ell)) \, . }
}\tag{4.9}
$$
\endroster
\endproclaim

\demo {\bf Proof of Proposition~4.5}
Taking small enough $\delta_k$ and $\theta_k$, the nonhyperbolic requirement of
(P.a) or (P.b) is obviously satisfied for the self-transition
on an $\omega$-petal by the  Hakim normal form \thetag{4.3}. 
We therefore concentrate on the hyperbolic condition for the
system \thetag{4.9}. Just like in [Ru2], the key is to reduce
to a Schwarz inclusion: 
\enddemo

\proclaim{Lemma 4.6 (Schwarz lemma contraction)} 
Proposition 4.5 holds, replacing condition (H)
for $(k,\ell) \in \SS^2 \setminus \SS_0^2$ by
$$
\eqalign{
&P^s_{\psi_j^{-1} (\xi_k)} 
(\hat f_{k\ell} ( R^\omega_k))\subset \Int (\DD^1_\ell) \, ,
 \quad P^u_{\psi_i^{-1} (\xi_\ell)} (\hat f^{-1}_{k\ell}( R^\omega_\ell))
\subset \Int (\DD^2_k) \, ,}\tag{4.10}
$$
and the hyperbolic condition for
$(k,k)\in \SS_0^2$ by the $P^u$-inclusion above in case (P.a) and
the $P^s$-inclusion in case (P.b). 
\endproclaim

\smallskip
The apparently weaker condition in
Lemma 4.6 implies (H): Indeed, the existence of a partial inverse 
$\phi_{k\ell,s }:\DD^1_k \times \DD^2_\ell
\to \Int (\DD^2_k)$, which is real-analytic in a neighbourhood
of $\DD^1_k \times \DD^2_\ell$, and is the unique
solution of
$
P^u_{\psi_j^{-1} (\xi_k)} 
\hat f_{k\ell} (w_1,\phi_{k\ell}(w_1,z_2))=z_2 
$,
can be obtained as in pp. 812-813 of [Ru2]. 
\qed
\medskip

The hard work consists now in proving Lemma 4.6, the technical 
but crucial  dynamical
lemma of this paper:

\demo{Proof of Lemma 4.6}
We first consider the case of a single  $F$-nonhyperbolic
fixed point $q$ of $f$ (which is $(0,0)$ in the charts), of multiplicity
$\nu+1=2$.

Let $\epsilon_0$ be small enough so that $V$, the $\epsilon_0$-neighbourhood
of the fixed point $q=0$, is contained in a chart of the atlas and in a domain
of definition of both Takens  and Hakim normal forms 
$f_H(x,y)$ and $f_T(s,t)$ from \thetag{4.2--4.3}. In particular, we skip
the chart index and do no distinguish between $\xi$ and
$\psi_j^{-1}(\xi)$ in the notation for this proof.
We shall  use the fact that for all $t\ge 0$
we have $\Re y(s,t) \ge (1-O(\epsilon_0) )t$ (the local strong stable manifold
is the same for both coordinates $(s,t)$ and $(x,y)$).
We also take $\epsilon_0$ small
enough so that $f^{\pm m}V\cap V\ne \emptyset$, for $m\ne 0$, is only
possible for very large
$|m|\ge m(\epsilon_0)$.

Pick a Markov partition $\QQ$ like described just after Definition~4.1,
ensuring that the rectangle containing $q$ is a subset of $V$.
For each $n \gg 1/ \epsilon_0$, we consider $\QQ_n$ the
$n$th refinement of $\QQ$ under $f$. We may replace $V$ by the
union of rectangles in $\QQ_n$ intersecting $V$, up to slightly
changing $\epsilon_0$.  We set $Q_q$ to be the rectangle containing $q=0$ in its (horizontal) boundary 
and  let $Q_{W^s(q)}$ be the set of rectangles of
$\QQ_n$ along the stable manifold of $q$, in particular
$Q_q \in Q_{W^s(q)}$. (Of course, if $n\gg m(\epsilon_0)$ then
$Q_{W^s(q)}$ winds back into $V$ and we get an infinite sequence of homoclinic
intersections.) We define $Q_{W^s_{loc}}(q)$ to be the union
of those rectangles in $Q_{W^s(q)}$ which are inside $V$.

Let us examine the rectangles $Q$ of $\QQ_n $ in $V$. 
Note that in the Takens normal form, the stable boundaries $\partial_s Q$
are horizontal segments, while the weak-unstable boundaries $\partial_u Q$
are curves which are close to vertical segments. 
We claim that the maximal diameter of $Q\in \QQ_n$ with $Q\subset V$
is $O(1/n)$, which is realised only for the vertical length of
rectangles of $V$ in $Q_{W^s(q)}$ (rectangles not along the
{\it global} stable manifold
have diameter $O(1/n^2)$).  This can be seen via the $f_T$
coordinate, since if $b^n(t_0)=\epsilon_0$ then $t_0=O(1/n)$, using 
the nonflatness Lemma~4.2
and $\nu+1=2$. In fact $O(1/n)$ is the diameter of $\QQ_n$
(also outside $V$).  Since $m(\epsilon_0)$ is large
and  we have $1-O(\epsilon_0)$ contraction outside of $V$,
the rectangles of $\QQ_n  \in Q_{W^s(q)}$ in  $V$
which are not along the {\it local} stable manifold
have (vertically realized) diameter at most $\epsilon(\epsilon_0)/n$
with $\epsilon\to 0$ as $\epsilon_0\to 0$.

\smallskip

Recall that $\gamma \le 1$ is the Hoelder smoothness of $\epsilon_2$.
If $\gamma\ge 1/2$ we may replace it by $0<\gamma' <1/2$, keeping the
notation $\gamma$.
Let $U$ be the union of elements of $\QQ_n$ which are
in an $n^{-\gamma}$ neighbourhood of $q=0$. 
We  next construct $\RR_n$ by modifying $\QQ_n$ in $U$.
Our aim is to ensure that we may choose a point $\xi_j$ in every rectangle 
$R_j$ in $\RR_n$, with $R_{j} \subset U$, and $q\notin R_j$ 
in such a way as to guarantee that if $R_k$ is another such rectangle, not along the local stable manifold of 
$q$, and $t_{jk}=1$, then
$f(\xi_j)=\xi_k$.

Let $\QQ_{n,0}$ be the set of rectangles 
$Q\in \QQ_n$ in $U$ such that $ f^{-1}(Q) \cap U=\emptyset$.
For each $Q_i$ in $\QQ_{n,0} \setminus Q_{W^s_{loc}(q)}$, consider all forward iterates  which
intersect $U$: $\{f^j (Q_i)\cap U\}$. 
For  $Q_i$ in $\QQ_{n,0}$ along the local stable manifold, we perform the same 
construction, except that
we set  $R'= f (Q_i)\cap Q_{W^s_{loc}(q)}$ and  $R= f (Q_i)\setminus Q_{W^s_{loc}(q)}$,
and we continue  iterating $R$, $R'$ until we leave $U$,  decomposing each
iterate which meets $Q_{W^s_{loc}(q)}$ into $R$ and $R'$.
The newly created sets $R_m$ are all Markov rectangles, and whenever $R_m \cap Q_i\ne \emptyset$ for
some $Q_i$ of the partition $\QQ_n$,
the complement $\widehat R=R_m \setminus( R_m \cap  Q_i)$ is also a Markov rectangle.
Letting $\QQ_{n,1}$ be the set of newly created complements $\widehat R$
such that $ f^{-1}(\widehat R) \cap U=\emptyset$,
we proceed as above, considering forward iterates in $U$ and taking appropriate intersections.
We repeat this procedure until $\QQ_{n,N}$ is empty.
Finally, we add to our collection $\{R_m\}$ of rectangles $\widehat R_q=Q_q\setminus \cup R_m$,
as well as $R_q= \overline {f(\widehat R_q) \setminus \widehat R_q}$,
and all its iterates $f^j(R_q)$ intersecting $U$. 

The rectangles $R_m$ are two by two disjoint and their union is $U$.
We define $\RR_n$ to be the union of the $R_m$s  and the rectangles of $\QQ_n$ outside
of $U$.
It is a Markov partition, which tends to be thinner horizontally and (slightly)
fatter vertically than $\QQ_n$  in $U$.
Note also (this is the announced feature) that we may choose a point $\xi_m=(s_\xi, t_\xi)$  in each
rectangle $R_m$ of $\QQ_{n,i}$, and a point $\xi_q$ in $R_q$, and consider the corresponding iterated points 
$\xi_\ell$ in elements $R_\ell$ of the partition $\RR_n$; when there
is an $R$, $R'$ bifurcation we ``follow'' the orbit in  $R$,
and take a new point in $\xi'\in R'$, making sure that $f(\xi) \in R$
and $\xi'$ are in the same $W^u$ leaf. If $R_k$ is not
adjacent to $W^s_{loc}(q)$ then $t_{\xi_k}\ge t_0=O(1/n)$ while if $R_m\ne \widehat R_q$
is adjacent to the local stable manifold of $q$, we take $\eta$ on the top $\partial_s$ boundary so that
$t_{\xi_m}\ge t_0=O(1/n)$ even in this case.

Note that inside $U$ we have
$|\Re x(s,t)| > (1-O( n^{-\gamma}) )|s| - O( n^{-\gamma})|t|$,
so that on the boundary of $U$ either $t= C n^{-\gamma}$ and thus
$\Re y \ge c n^{-\gamma}$ or $t< C n^{-\gamma}$ and $|\Re x |\ge  C|s|$
with $|s|\ge cn^{-\gamma}$. Hence, by \thetag{4.5}, the contraction factor outside of
$U$ (and also when
entering $U$ from outside of $U$, or just when leaving $U$)
is $1- c n^{-\gamma}$ (note that the points involved belong to
the rectangles and are thus real points).
Replacing $U$ by $V$, the same argument gives a contraction factor $1-C \epsilon_0$.

If $t_{jk}=1$ and both rectangles $R_k$ and $R_j$ are
outside of $V$
we can   apply the construction in  the lemma of [Ru2, p. 811].
Let us recall here the key estimate involved:
Fixing $\xi\in R_j$ and $\eta\in R_k$ we first observe (see [Ru2, p. ~813])
that $|f^{-1}(\eta)-\xi|_\xi \le O(1/n)$.
Thus, using \thetag{4.4} and \thetag{4.6, 4.8} 
(in particular  $\epsilon^F_1(a)=O(a)$),
$$
\eqalign{
\| f^{-1}(v )-\xi \|^u_\xi &\le  \| f^{-1}(\eta)-\xi\|^u_\xi+
\epsilon_2(d(\xi,f^{-1}(\eta))\cdot  |f^{-1} (v)- f^{-1}(\eta) |_{f^{-1}(\eta)}\cr
&\qquad+
 \|f^{-1}(v) -f^{-1}(\eta) -Df^{-1}_\eta(v-\eta)\|^u_{f^{-1}(\eta)}+ 
  \|Df^{-1}_\eta (v-\eta)\|^u_{f^{-1}(\eta)} \cr
&\le C | f^{-1}(\eta)-\xi|_\xi+
C  |f^{-1} (\eta)-\xi|^\gamma_{f^{-1}(\eta)} |f^{-1}(v)-f^{-1}(\eta)|_{f^{-1}(\eta)} \cr
&\qquad+ 
  \epsilon^F_1(|v-\eta|_{\eta}) |v-\eta|_{\eta} +  (1-C_\nu \epsilon_0) \|v-\eta\|^u_\eta \cr
&\le  {C \over n} +
{C \over n^\gamma} |v-\eta|_\eta 
+  \widetilde C |v-\eta|_\eta^2   + (1-C_\nu \epsilon_0) \|v-\eta\|^u_\eta  \, . \cr
}
$$
Taking  the size $\delta_k=\delta_j$
of the $\omega$-rectangles
$R^\omega_{j,k}$  to be $C \epsilon(\epsilon_0) n^{-\gamma}$ (this choice  will turn out to be useful
later), we get
$
\|v-\eta\|_\eta\le \delta_k \Longrightarrow 
\|f^{-1}(v)-\xi\|^u_\xi< \delta_j 
$.
The $P^s$-inclusion is similar (in fact easier), we shall concentrate
on the $P^u$-inclusion.

If $t_{jk}=1$ and both rectangles $R_k$ and $R_j$ are in $V$, but
outside of $U\cup Q_{W^s_{loc}(q)}$, the $\|\cdot\|^u$- contraction
of $f^{-1}|_{R_k \cap f(R_j)}$ is at least $1-C_\nu n^{-\gamma}$.
We can essentially apply the above estimate, using also that
 $|f^{-1}(\eta)-\xi|_\xi \le \epsilon/n$ in this case:
$$
\eqalign{
\| f^{-1}(v )-\xi \|^u_\xi 
&\le C | f^{-1}(\eta)-\xi|_\xi+
C  |f^{-1} (\eta)-\xi|^\gamma_{f^{-1}(\eta)} |f^{-1}(v)-f^{-1}(\eta)|_{f^{-1}(\eta)} \cr
&\qquad+ 
C  \epsilon_1^F(|v-\eta|_{\eta}) |v-\eta|_{\eta} +  
(1-C_\nu n^{-\gamma}) \|v-\eta\|^u_\eta \cr
&\le { C \epsilon\over n} +
\left ( {C \epsilon \over n} \right )^\gamma  |v-\eta|_\eta 
+  \widetilde C |v-\eta|_\eta^2   + (1-{C_\nu \over n^{\gamma}}  )\|v-\eta\|^u_\eta  \, , \cr
}\tag{4.11}
$$
(note that $\eta$ and $f^{-1}(\eta)$ are in $\Lambda
\subset \widetilde \Lambda$). Taking  the size $\delta_k=\delta_j=\delta$ to be 
$\delta= C\epsilon(\epsilon_0) n^{-\gamma}=\hat \epsilon n^{-\gamma}$, we get
$
\|v-\eta\|_\eta\le \delta_k \Longrightarrow 
\|f^{-1}(v)-\xi\|^u_\xi< \delta_j 
$.

\medskip

In order to obtain the  $P^u$-inclusion for $R_j$ or $R_k$ in $U$, 
note first that any $R_m$ in $U\setminus Q_{W^s_{loc}(q)}$
lies entirely between two horizontal lines
$t\equiv b^\ell(t_0)$ and  $t\equiv b^{\ell+1}(t_0)$
for an integer $\ell=\ell(R_m)$ between $0$ and $n-1$
(recall the definition of $t_0=O(1/n)$).
Using Lemma~4.2, it is
not very difficult to see that the vertical diameter of  such
a rectangle $R_m$  is 
not larger than $B(n)(b^{\ell( R_m)}(t_0))^2$,
for some constant $B(n)\ge 1$ tending to $1$ when $n\to \infty$
(for fixed $\epsilon_0$).
Indeed, recall $n\gg 1/\epsilon_0$ 
and notice that if $b^{\ell}(t_0)\le n^{-\gamma}$ then $\ell\le\ell_{max} \le  n - A n^\gamma$
for some  $A > 0$.  For such $\ell$ we have
$$
0< b^{(\ell+1)}(t_0)-b^{(\ell)}(t_0)=
B(n)  (b^{(\ell)}(t_0))^2 \le \biggl (1+O\left ( {1\over n^{\gamma}}\right )\biggr )
 (b^{(\ell)}(t_0))^2\, .
\tag{4.12}
$$
Along the local stable manifold of $0$,   the diameter is of course  $t_0=O(1/n)$.

\smallskip

Discuss first rectangles $R_j$, $R_k$ in $U$
with $t_{jk}=1$, such that $R_j$ is  {\it not} adjacent to the
local stable manifold of $q$.
Take as reference points for the (complex) $\omega$-rectangles the
chosen points $\xi=\xi_j$ and $\eta=\xi_k$.  Note that
$\ell(R_j)=\ell(R_k)-1\ge 0$ and set $\ell_{k/j}=\ell(R_{k/j})$. 
Then, since $f^{-1}( \eta)=\xi$,  we only have two nonzero terms out of 
four in \thetag{4.11}, and if $\|v-\eta\|_\eta \le \delta_k$,
$$
\eqalign{
\| f^{-1}(v )-\xi \|^u_\xi &= \|f^{-1}(v) -f^{-1}( \eta)\|^u_{f^{-1}(\eta)}\cr
   &\le \|f^{-1}(v) -f^{-1}(\eta) -Df^{-1}_\eta(v-\eta)\|^u_{f^{-1}(\eta)}+ 
       \|Df^{-1}_\eta (v-\eta)\|^u_{f^{-1}(\eta)}\cr
   &\le \widetilde  C |v-\eta|_\eta^2   + (1- C_\nu b^{\ell_j}(t_0)) \|v-\eta\|^u_\eta \cr
&\le (\widetilde  C'  \delta_k + 1- C_\nu b^{\ell_j}(t_0)) \delta_k \, .  \cr}
\tag{4.13}
$$

Since $C'$ does not depend on $\epsilon_0$ or $n$, and
$C_\nu > \nu\ge 1$ while $B(n)\to 1$ when $n\to \infty$,
for each $\epsilon_0$, up to taking smaller
$\epsilon_0$, we may assume that $\widetilde C' \hat \epsilon(\epsilon_0)(1+B(n))-C_\nu< -B(n)$.
We take 
$$
\delta_m= \hat \epsilon b^{\ell_m}(t_0)\, , \quad m=k\, , j\, .
$$

Also, for $\ell_k=0, \ldots, \ell_{max}$ we may assume $\delta_k\le \hat \epsilon  n^{-\gamma}$.
Then, on the one hand  the real projection of $R^\omega_{k/j}$  contains  $R_{k/j}$,
and on the other, if $\|v-\eta\|_\eta \le \delta_k$,  then by \thetag{4.13}
$$
\eqalign{
\| f^{-1}(v )-\xi \|^u_\xi  
& \le
\bigl (1 +( \widetilde C'\hat \epsilon(1+B(n) b^{\ell_j}) -C_\nu) b^{\ell_j}(t_0)\bigr ) 
\cdot \hat \epsilon b^{\ell_k}(t_0)\cr
& \le
\bigl (1-B(n) b^{\ell_j}(t_0) \bigr )\cdot
\bigl (1+B(n)b^{\ell_j}(t_0) \bigr )  \cdot \hat \epsilon b^{\ell_j}(t_0)
< \delta_j  \, ,
}
$$
so that we have the required Schwarz lemma inclusion property.
If $R_j$, $R_k$ in $U$
with $t_{jk}=1$, are such that $R_k$ is  {\it not} adjacent to the
local stable manifold of $q$, while $R_j \ne \widehat R_q$ is adjacent to it,
we have $\ell_k=0$, so that we already set $\delta_k=\hat \epsilon t_0$.
We take $\delta_j =C_j /n$, where $C_j$ is a constant (independent
of $n$) to be made more precise later on. So,  if $\|v-\eta\|_\eta \le \delta_k$, 
$$
\eqalign{
\| f^{-1}(v )-\xi \|^u_\xi  &\le
\bigl (1 +( C'\hat \epsilon-C_\nu) t_0\bigr ) \cdot \hat \epsilon t_0\le C_j/n  \, .
}
$$
If $R_j=\widehat R_q$, we have $\ell_k=0$ and $\delta_k=\hat \epsilon t_0$
and we need to show that $\|v-\eta\|_\eta \le \delta_k$ implies that
$f^{-1}(v)$ belongs to an $\omega$-petal at $q$. Taking the size of this
$\omega$-petal to be $C_q/n$ for $C_q \ge 1$ but not very large, this follows from the
fact that the $t$- (and thus the real part
of the $y$-) coordinate of $f^{-1}(v)$ is positive, but not much bigger than
$t_0$.

For the transitions $t_{jk}$ from $R_k$, outside of $U$, to $R_j$, inside of $U$, we have
as mentioned above a contraction  factor $1-O( n^{-\gamma})$,
and we may use the ``ordinary" four-term estimate \thetag{4.11} 
since $\delta_j= \hat \epsilon b^{\ell_{j_{max}}}(t_0)$
while $\delta_k= \hat \epsilon n^{-\gamma}$. 
For the transitions from $U$ (not along the local
stable manifold) to outside of $U$, we apply
\thetag{4.5} as in \thetag{4.11}, using $\delta_{k} \le \hat\epsilon n^{-\gamma}$,
and the contraction factor $1-C_\nu n^{-\gamma}$.

We finally discuss transitions in $V$ along the local stable manifold of $q$.
Let $R_j$ be a rectangle with $f(R_j)\subset \widehat R_q$. 
If  $\|v-q\|_q\le 2\delta_q=2C_q/n$ (this is in fact a weaker condition
than the petal condition) we have from the usual four-term
estimate, but without any contraction factor,
$$
\| f^{-1}(v) -\xi \|^u_\xi\le  O(1/n)+ O(1/n^{\gamma}) |v-q|_q+
  C |v-q|^2_q   +  \|v-q\|^u_q
 \le  C'_q/n \, .\tag{4.14}
$$
We take $\delta^u_j=C'_q/n$, which defines
the constant $C_j$ for this rectangle. 
When we iterate, progressing to the
left or to the right along  of $W^s_{loc}(q)$, we get (recall that $\xi$ and
$f^{-1}\eta$ are on the same unstable leaf, and also that $t_\xi=t_0$ while
$b(t_{f^{-1}(\eta)})=t_\xi$)
$$
\| f^{-1}(v) -\xi \|^u_\xi\le  O(1/n^2)+ O(1/n) |v-\eta|_\eta+
  C |v-\eta|^2_\eta   + (1-c/n) \|v-\eta\|^u_\eta
 \le(1 +O(1/n)) \delta_k\, , \tag{4.15}
$$
and the diameter $\delta_j=(1 +O(1/n))\delta_k$ grows.
However, the number $n_0$ of iterations of $f^{-1}$ to escape
from $V$ to the left, starting from a component of
$f^{-1}(\widehat R_q)\setminus R_q$  is at most  $n$. So the cumulated
factor is smaller than $\zeta_n=(1+C/n)^{n}$, 
which is uniformly bounded as $n \to \infty$, by $\zeta\gg 1$, say. 
We may choose all the $C_\ell \le \zeta C'_q$.
Finally, 
since $\zeta O(1/n)\le \hat \epsilon  n^{-\gamma}$, we get
the required inclusion when we leave $U$, and, a fortiori, $V$. 

\smallskip

If there are several $F$-nonhyperbolic fixed points (of the same
multiplicity $\nu+1=2$), the construction
above works too.  The $E$-nonhyperbolic case is analogous.
Dealing with
coexistence of $E$ and $F$-nonhyperbolic points 
(of the same multiplicity) does not cause any 
problems, since we never used the exponential smallness
in the strong direction in our estimates: $O(1/n)$ was enough outside of $V$
and on $Q_{W_{loc}(q)}$
and $O(\epsilon/n)$ elsewhere.
\qed\enddemo

\remark{Remark 4.7} When $n \to \infty$, the cardinality
of $\SS_1$ goes to infinity but the cardinality of $\SS_0$ is fixed.
Also, each  map $\hat f_{kk,n}$ associated to a
nonhyperbolic fixed point (i.e., $k \in \SS_0$) for $n\ge 1$
is just the restriction of $\hat f_{kk,1}$ to  smaller
domains $\DD^1_k$, $\DD^2_k$, which can be constructed
in a way compatible with the definition of $\BB_0(\epsilon)$ in
Lemma~3.1\therosteritem{2}.
\endremark

\head 5. Spectral harvest 
\endhead

\subhead 5.A Proof of Theorem A in the parabolic case
\endsubhead

In this section, we put together  the results of Sections 2--3 and Section 4
to prove Theorem A, under the additional
assumptions that the nonhyperbolic periodic points which are not
$\Omega\setminus P$ isolated are fixed, 
never have an eigenvalue equal to $-1$, and always have
multiplicity $\nu+1=2$.

As explained in Section 4, we may restrict to a single basic set
$\Lambda_j$.  Let $\NN$ be the set of nonhyperbolic non $\Omega\setminus P$ isolated periodic
orbits and recall the definition \thetag{1.3} of  $\Sigma(p)$ 
for $p\in \NN$. We fix an arbitrary open neighbourhood  $\WW$ of
$\cup_{p\in \NN \cap \Lambda_j} \{ 1/z \mid z \in \Sigma(p)\}$, and we
shall show that $d_{f|\Lambda_j}(z)$ is holomorphic in
$\complex \setminus \WW$.

Putting together Proposition~4.5 and Lemmas 3.8 and 3.1--3.2, with
Remark 4.7,
and using Manning's counting trick [M] (see also [Ru2, p. 817]), 
we obtain that for each  Markov partition $\RR_n$ of small
enough diameter,
the determinant may be written as a ratio of two  functions 
which are analytic  in $\complex \setminus \WW$:
 $$
d_{f|\Lambda_j}(z) = { d_1(z) \over d_2(z)} \, .
$$
Indeed, $d_1$ is the regularised Fredholm determinant of the transfer
operator of the hyperbolic analytic map $\hat f$
associated by Proposition~4.5 to a Markov partition $\RR_n$, while $d_2$
is the determinant of a second symbolic map
$\hat f_2$, associated to the auxiliary subshift corresponding
to the pairs of adjacent rectangles in $\RR_n$.
(Note that the nonhyperbolic periodic orbits are not counted
at all and therefore not overcounted, while the hyperbolic periodic
orbits on boundaries of rectangles of the partitions coinciding with
boundaries of our basic set are
only counted once.)

It thus only remains to check that all zeroes of $d_2(z)$ are
cancelled by zeroes of $d_1(z)$, for each $n$. 
The argument for this uses our assumption 
that $M$ is two-dimensional:

\proclaim{Sublemma 5.1 (Disjoint boundary periodic orbits)}
Let $\Lambda$ be a basic set of  a $C^2$ diffeomorphism $f$ on a compact
surface having dominated splitting over its nonwandering set. Then for 
each $\epsilon > 0$ there are two finite Markov
partitions of diameter smaller than $\epsilon$, and  such that
the sets $B_1$, $B_2$ of periodic orbits of $f$
lying on their respective boundaries may only intersect on the
boundary of $\Lambda$.
\endproclaim

See the Proposition on p.~817 of [Ru2],
to which we refer for a proof of Sublemma ~5.1, valid in the Anosov case
(there, $B_1$ and $B_2$ can be taken disjoint), and which may
be extended to Axiom A and  dominated splitting.

\demo{End of the proof of Theorem A (parabolic case)}
Sublemma 5.1 finishes the proof. Indeed
we may take two sequences $\RR_{1,n}$, $\RR_{2,n}$ of Markov
partitions so that, on the one hand,  $\GG(\epsilon)$ in
Lemma 3.2 goes to zero, and on the other, the
sets $B_1^n$ and $B_2^n$ satisfy the properties in Sublemma~5.1.
The argument on p. 818 of [Ru2]  applies for each  $n$.
$\qed$\enddemo

\subhead 5.B The general case
\endsubhead

We discuss briefly the changes needed to handle  the general
case. First note that
the arguments in the appendix apply to  neutral eigenvalues
$-1$, working on
both sides of the central manifold, and to periods larger than
one, by considering the corresponding iterate of $f$. 
Let us now explain the changes in Sections 2--3 and 4.

\smallskip

Regarding $\lambda_{E,F}=-1$, we may easily generalise
(P.a), (P.b) by allowing normal forms $-z_2 -z_2^{1+\nu}
+ h.o.t.$ instead of \thetag{2.3} and $-w_1-w_1^{1+\nu}+h.o.t.$
instead of \thetag{2.4}. We then use two petals, one on each side of $0$.
The square of the
local  (P) dynamics has neutral eigenvalue $+1$ and the same multiplicity
$\nu+1$.
We may thus adapt the proof of Lemma~ 3.1, using the fact that
the spectrum of $(\TT^{\pm})^2$ is $[0,1]$ so that the spectrum of
$\TT^\pm$ is contained in $[-1, 1]$. We get the sets announced in \thetag{1.3}
for $P=1$. In Section~4, we have to take into account the fact that
the dynamics oscillates between both sides of the strong local manifold,
but this does not require any serious changes.

\smallskip

Regarding {\it fixed points of  multiplicities $\nu+1\ge 3$,}  using the conformal change
of variables $(x, y)\mapsto (x, y ^\nu)$ as in [Ru3, \S 2.6],
we may adapt the contents of Sections 2 and 3.
(Note that for even $\nu$ we need to consider two real petals for
each $j\in \SS_0$, while for odd $\nu$ we have just one petal.)
We must also adapt the proof of Lemma 4.6.
The main changes are the weaker contraction $1-C_\nu t^\nu$
(note however that $\epsilon_1^\nu$ has stronger decay, see
\thetag{4.7}), and
the larger diameter of the refined partition: $O(1/n^{1/\nu})$ in
general, and  $O(\epsilon(\epsilon_0)/n^{1/\nu})$ (except along the
local strong manifold) in an $O(\epsilon)$ neighbourhood $V$ of
the nonhyperbolic periodic point $p=0$. To deal with this,
we replace $\gamma$ by
$\gamma'<1/(\nu(\nu+1))$ if $\gamma \ge 1/(\nu(\nu+1))$.
Then, $\nu \gamma+ \gamma < 1/\nu$.
We take $U$ to be an $O(n^{-\gamma})$ neighbourhood of
$0$ (for the new, possibly smaller $\gamma$) and choose $\delta_m=C\epsilon/n^\gamma$
outside of $U$. Then, \thetag{4.11} becomes
$$
\eqalign{
\| f^{-1}(v )-\xi \|^u_\xi 
&\le { C \epsilon\over n^{1/\nu}} +
\left ( {C \epsilon \over n^{1/\nu}} \right )^\gamma  |v-\eta|_\eta 
+ { \widetilde C\over n^{(\nu-1)\gamma}} |v-\eta|_\eta^2   + 
(1-{C_\nu \over n^{\nu\gamma}}  )\|v-\eta\|^u_\eta  \, , \cr
}\tag{5.1}
$$
so that $\|v-\eta\|_\eta\le \delta_k$ implies $\|f^{-1}(v)-\xi\|^u_\xi\le \delta_j$.
Then,  \thetag{4.12} becomes
$$
\eqalign{
0< b^{(\ell+1)}(t_0)-b^{(\ell)}(t_0)
&\le B(n) (b^{(\ell)}(t_0))^{\nu+1}  \le \bigl (1+O( n^{-\gamma})\bigr )
 (b^{(\ell)}(t_0))^{\nu+1} \, .}
\tag{5.2}
$$
For the the two-term transitions from $R_k$ to $R_j$ within $U$, we
set $\delta_m=\hat \epsilon b^{\ell_m}(t_0)$, and  invoke \thetag{4.7} inside the
$b^{(\ell_{k}+1)}(t_0)$ neighbourhood. Then 
we may assume that $\hat \epsilon \widetilde  C' ( 1+B(n))^{2\nu-1}    - C_\nu
< -B(n)$ and \thetag{4.13} becomes
$$
\eqalign{
\| f^{-1}(v )&-\xi \|^u_\xi 
   \le \widetilde  C  ( b^{\ell_k+1}(t_0))^{\nu-1}  |v-\eta|_\eta^2   
+ (1- C_\nu (b^{\ell_j}(t_0)) ^\nu ) \|v-\eta\|^u_\eta \cr
&\le(1+ (\widetilde  C' ( b^{\ell_k+1}(t_0))^{\nu-1}  \delta_k 
- C_\nu (b^{\ell_j}(t_0))^\nu) \delta_k\cr
&\le 
[1+
(\hat \epsilon \widetilde  C' ( 1+B(n))^{2\nu-1}   - C_\nu )
(b^{\ell_j}(t_0))^{\nu}
]
(1+B(n) ( b^{\ell_j}(t_0))^\nu)  \delta_j < \delta_j\, .  \cr}
\tag{5.3}
$$
Along the local strong manifold, we take $\delta_m=C_m/n^{1/\nu}$,
and the diameter grows 
$$
\eqalign{
&\| f^{-1}(v) -\xi \|^u_\xi\cr
&\qquad\le  O(1/n^{1+1/\nu})+ O(1/n^{1/\nu}) |v-\eta|_\eta+
  C |v-\eta|^2_\eta /n^{1-1/\nu}  + (1-C_\nu/n) \|v-\eta\|^u_\eta\cr
 &\qquad\le(1 +O(1/n^{1/\nu})) \delta_k\, , }\tag{5.4}
$$
but (up to considering as a single Markov rectangle at $p$ all the
rectangles to the left and to the right in a horizontal neighbourhood
of horizontal size $\epsilon /n^{\nu}$, which does not interfere
with the other computations) we may assume that the number of iterations $n_0$ within
$U$ is smaller than
$C\log n \ll n^{\gamma/\nu}$.
Since $\zeta  /n^{1/\nu} \ll 1/n^\gamma$, we are done.

Finally, if there are nonhyperbolic fixed points of different indices,
we set $\nu=\max \nu_j$ and use the same $\gamma' \le \gamma$
with $\gamma'< 1/(\nu(\nu+1))$ for all neighbourhoods
$V_j$ and $U_j$.

\smallskip

If we have periodic points of {\it periods $P\ge 2$ in $\NN$,} they will be associated
to a periodic cycle in the symbolic set $\SS_0$ (work first with $f^P$). The corresponding
symbolic transfer operator $\widehat \LL_0$ will have associated
periodic blocks of periods $P_i$. Its spectrum is thus contained in the
union over $i$ of the $P_i$th roots of the spectrum of the block $\widehat \LL_{0,i}^{P_i}$.
This gives the $P$-th roots announced in \thetag{1.3} and thus the slits
in Theorem ~A.\qed

\smallskip

\head Appendix: Constructing  an adapted metric
\endhead

Let $f$ be a real-analytic diffeomorphism of a compact surface
having a dominated splitting $T_\Omega M=E\oplus F$.
Let $\Lambda$ be a basic set from the decomposition
[PS2] of the nonwandering set $\Omega$.
We consider an analytic atlas  for $f$, as in Section 4
and denote by $|\cdot |_z$ the riemannian norm induced on complex charts.
(We systematically drop the chart index in this appendix.)
Let $\NN$ be the finite set of nonhyperbolic periodic points of $f$
in $\Lambda$. For each $q_j\in \NN$, we denote by $\nu_j+1\ge 2$ 
its index, i.e., the
order of the zero $f-\Id$ at $q_j$ in the charts
(in other words, the multiplicity of $f-\Id$, recall from Section 4 that this
multiplicity is finite). In this appendix, we allow all $\nu_j\ge 1$, however, for simplicity, 
we assume that all
points in $\NN$ are {\it fixed points} and that their nonhyperbolic multiplier
is equal to $+1$ (not $-1$). See Section~5.B for the general case.

Our aim in this appendix is to prove Lemma 4.3, i.e.,
to construct two adapted semi-norms
(see [Ru2 pp. 809--810] for a brief  account of the
hyperbolic case)   $\|\cdot \|_{E,z}=\|\cdot \|^s_{z}$, 
$\| \cdot\|_{F,z}=\| \cdot\|^u_{z} $
on the complex tangent bundle $T_{\Lambda^\complex}M$ over
a complex neighbourhood $\Lambda^\complex$ of $\Lambda$, and
an admissible complex extension $\widetilde \Lambda$ 
of $\Lambda$, such that \thetag{4.4--4.8} hold.

We shall recycle some ideas of Crovisier [Cr, \S 5.3]
(introducing simplifications arising from [PS1, PS2]),
but we must modify his
construction which uses smooth (not analytic) coordinates 
(see [Cr, \S4.2]), since we need to control the complex extensions.
Another difference is that Crovisier constructs an
adapted metric in the sense that 
$
\|Df^{-1}_{/F(z)}\|^C_z < 1 $ and 
$\|Df _{/E(z)}\|^C_z<1
$ if $z \notin \NN$, while we need the quantitative
estimate \thetag{4.5}. 

\medskip 

{\bf Adapted metrics close to $\NN$  --- The first global norm $\|\cdot \|'_{z}$}

Our first step is to construct  local semi-norms $\|\cdot\|'_{F/E,z}$  in a
neighbourhood of each nonhyperbolic fixed point  of
$f|_\Lambda$. 

Recall (see e.g. [PS1, Lemma 3.2.1]) that for suitable $\alpha <1$ the (complex) cone
$$
\CC^F_{\alpha,z}=\{ w \in T^\complex_z \Lambda \mid
w = u + v \, ,  u \in E(z) \, , v \in F(z) \, ,
|u| \le \alpha |v|\}\, \tag{A.1}
$$ 
is invariant under $Df$, for $z$ in a complex neighbourhood 
of $\Lambda$.

\proclaim{Lemma A.1 (Local semi-norm)}
Let $f$ be a real-analytic surface diffeomorphism with dominated
splitting over $\Omega$,  and $\Lambda$ 
a basic set.
Let $q\in \Lambda$ be  an $F$-nonhyperbolic fixed point of index $\nu+1\ge2$.
Then there exist a semi-norm $\|\cdot\|'_{F,z}$ over a complex neighbourhood 
$B^\complex$ of
a (real) neighbourhood  $B$ of $q=0$, 
and an admissible complex extension $\widetilde B$ of $B$,
so that \thetag{4.5} holds for all $w\in \CC_{\alpha,z}^F$,
i.e., for some $C_\nu >  \nu$ 
 $$
\|Df^{-1}w\|'_{F, z}\le (1-C|z_E|)(1-C_\nu |\Re (z_F^\nu) |)) \|w\|'_{F,f(z)}\, , \, \forall
z \in  \widetilde B
\, ,  \forall w\in \CC_{\alpha,z}^F \, .\tag{A.2}
$$ 
Furthermore $z\mapsto \|\cdot \|'_{F,z}$ is $C^\infty$,
\thetag{4.7} holds on $B^\complex$ if $(z-z_0)\in \CC_{\alpha, z_0}$, and
\thetag{4.8} holds. We have analogous statements
if $q$ is $E$-nonhyperbolic.
\endproclaim

\demo{Proof of Lemma A.1}
As mentioned in the proof of Lemma ~4.2, by [H], 
we may express $f$ in real-analytic coordinates in a neighbourhood of $0$ as
$$
f_H(x,y)=(g(x) +yh(x,y)\, , y+y^{\nu+1}+Ay^{2\nu+1}+
\sum_{j=2\nu+2}^\infty a_{j}(x)y^{j}) \, , 
$$
where $A$ is a    constant, and
$g$, $a_j$, and  $h$, respectively, are real-analytic in a neighbourhood of
$0$ in $\complex$, respectively  $\complex^2$, with  $g(0)=0$, $h(0,0)=0$, and $0<|g'(0)|<1$.
Writing $h_x=\partial_x h$,  $h_y=\partial_y h$, we get 
$$
\eqalign{
Df_{H,(x,y)}(u,v)&=
((g'(x)+yh_x(x,y))u+(h(x,y)+ y h_y(x,y))v,\cr
&\qquad\qquad\qquad\qquad\quad 
\sum_{j= 2\nu+2}^\infty (a_{j}'(x)y^{j}) u +(1+(\nu+1)y^\nu+\cdots) v)\, .}\tag{A.3}
$$ 

Note  that   if $\Re (y^\nu) \ge 0$
$$
\eqalign{
|(1+(\nu+1)y^\nu+O(y^{2\nu}))  v|\ge
 (1+  C_\nu |\Re (y^\nu)| )  |v|\, . }\tag{A.4}
$$

We may view $0$ as a partially hyperbolic fixed point of a
$\CC^\infty$ map on a neighbourhood of $0$ in $\real^4$  and use the
Takens [T] standard coordinates $(\vec s, \vec t)=(s_1, s_2, t_1,t_2)$. Here, 
$f_T(\vec s, \vec t)=(G_{\vec t} (\vec s), H(\vec t))$ where 
$G_{\vec t}$ is a contraction, i.e., both its eigenvalues have absolute
values $<1$, uniformly in $t$, while the two eigenvalues of $H$ have moduli $1$.
In particular (up to taking
a smaller neighbourhood) there is $\gamma <1$ so that
$
\sup_{\vec t} |G_{\vec t}(\vec s)|\le \gamma |\vec s|
$.
For small  $\beta >0$,
define a semi-norm for a complex vector $(u,v)$ over $z=(x,y)$ by
$$
\|(u,v)\|'_{F,(x,y)}= (1-\beta |\vec s(x,y)|)| v| \, .\tag{A.5}
$$

Take $c < \beta (1-\gamma)$.
By definition, and by \thetag{A.3--A.4}, we have  for $\Re y^\nu \ge 0$  and $(u,v)\in \CC^F_\alpha$
$$
\eqalign
{
&\|Df_{H,(x,y)}(u,v)\|'_{F,f(x,y)}\cr 
&\qquad\qquad=(1-\beta|\vec s(f_H(x,y))|) \cdot | \sum_{j\ge 2\nu+2}(a_{j}'(x)y^{j}) u 
+(1+(\nu+1)y^\nu+\cdots) v|\cr
&\qquad\qquad\ge  (1-\gamma \beta |\vec s(x,y)| )
(1+ C_\nu |\Re (y^\nu)|) |v|   \cr
&\qquad\qquad\ge (1+c |\vec s(x,y)|)  (1+C_\nu |\Re (y^\nu)|)\|(u,v)\|'_{F,(x,y)}\cr
&\qquad\qquad\ge (1+\tilde c|x|) (1+ C_\nu|\Re (y^\nu)|) \|(u,v)\|'_{F,(x,y)}\, .\qed
}
\tag{A.6}
$$
\enddemo

Define a local norm in a neighbourhood of an $F$-nonhyperbolic
fixed point $q\in \Lambda$:
$$
\|(u,v) \|_z':=\cases 
\|(u,v) \|'_{F,z} &\hbox{if } (u,v)\in \CC^F_\alpha\, , \cr
|u|=: \|(u,v) \|'_{E,z} &\hbox {otherwise}\, .
\endcases 
$$

It is not difficult to check, using Lemma~A.1 and the cone property,
that the new norm is adapted in $\widetilde B$ in the sense that
$$
\cases
\|Df _{/E(z)}\|'_{z} < 1 \, , \forall z \in  B^\complex \, , \cr
\|Df^{-1}_{ /F(z)}\|'_{z} \le 1 \, , \forall z \in \widetilde B \, ,
\hbox{ with equality if and only if } z =q\, . \cr
\endcases
\tag{A.7}
$$

We may  extend the norm $\|\cdot \|'_{z}$ continuously  to $\widetilde \Lambda$ by glueing
it with the riemannian norm $|\cdot|_z$.
We replace the riemannian norm by this new equivalent norm $\|\cdot\|'_z$. Of course, 
we still have a dominated splitting.

\smallskip

{\bf Controlling hyperbolicity away from $\NN$ --- The second global norm $\|\cdot \|''_{z}$}

We will  use the following results proved in
[PS1] and [PS2]:

\proclaim{Proposition A.2 [PS1, Theorem B]}
Let $f$ be a $C^2$-diffeomorphism on a compact
surface, and $\Lambda$ a compact $f$-invariant set having a dominated
splitting. Assume that all the periodic points in $\Lambda$ are
hyperbolic of saddle-type and that $\Lambda$ does not contain normally invariant
curves.  Then, $\Lambda$ is hyperbolic.
\endproclaim

\proclaim{Lemma A.3 [PS2, see Proposition 3.1]}
Let $f$ be a $C^2$-diffeomorphism on a compact
surface, and $\Lambda$ a compact $f$-invariant set having a dominated
splitting. Let $q$ be a periodic point of $\Lambda$ and let $z_0\in
\Lambda$ be a point in the local unstable manifold of $q$. Then,
there are a neighbourhood $V$ of $z_0$ and an integer $k_0$ such that if $z\in
V\cap\Lambda$  and $k> k_0$ then 
$\|Df^{-k}_{/F(z)}\|' <\frac{1}{2}$.
\endproclaim

We will also use the following two lemmas (and their  versions
exchanging $E$ and $F$):

\proclaim{Lemma A.4} Let $f$ and $\Lambda$ be as in Lemma~A.3 and 
assume that $q\in \Lambda$ is an  $F$-nonhyperbolic fixed point,
and is the only nonhyperbolic periodic point of $f$ in $\Lambda$. 
Then, there is
$m_0\ge 1$ such that for any $z\in\Lambda\setminus\{q\}$ 
there is $1\le n(z)\le m_0$ such that $\|Df^{-n(z)}_{/F(z)}\|'< 1$.
We shall take $n(z)$ minimal with this property.
\endproclaim

\demo{Proof of Lemma A.4} 
If $z$ is in some (real) neighbourhood $B$ of $q$, we may take
$n(z)=1$ by the construction after Lemma~A.1.

Assume for a contradiction that the  conclusion of the lemma does not 
hold (outside of $B$). Then there
exists a sequence $z_n \notin B$ such that $\|Df^{-j} _{/F(z_n)}\|'\ge 1$,
for all $1\le j\le n$. Taking $z$ to be an accumulation point of the $z_n$, we have
$z\notin B$ and $\|Df^{-j}_{/F(z)} \|'\ge 1$ for all $j\ge 1$.
Now, if $q\notin\alpha(z)$ then $\alpha(z)$ is a compact set with
all  periodic points  hyperbolic, and from
Proposition ~A.2  this set must be hyperbolic,  a contradiction.

Consider now the case  $q\in \alpha(z)$. We
analyse two situations:  $\{q\}=\alpha(z)$ or
$\{q\}\varsubsetneq\alpha(z)$. In the first case, 
$z\in W^u(q)$ and thus there is $k_0$ such that $f^{-k_0}(z) \in
W^u_{loc}(q)$. Therefore
$\|Df^{-j}_{/F(f^{-k_0}(z_))}\|'\to 0$, and so also
$\|Df^{-j}_{/F(z)}\|'\to 0$. Then by continuity,  taking $z_n$
close enough to $z$, we also get a contradiction. Finally, if
$\{q\} \varsubsetneq\alpha(z)$, we may take a neighbourhood $V$ in the
local unstable manifold of $q$ and a subsequence
$f^{-j_n}(z_n)\in V$ with $j_n$ uniformly bounded.  
From Lemma~A.3, we also get a contradiction. 
\qed
\enddemo

\proclaim{Lemma A.5}  Let $f$, $\Lambda$ be as in Lemma~A.4. 
For each neighbourhood $B_1$ of $q$
and every $0<\lambda_1 < 1$
there exists $m_1$ and a complex neighbourhood $W_1$ of
$ \Lambda\setminus B_1$ such that for each
$z\in W_1$ then $\|Df^{-m_1}_{/F(z)}\|'<\lambda_1$.
\endproclaim

\demo{Proof of Lemma A.5}
For arbitrary $z\in \Lambda\setminus \{ q\}$ let  $n(z)\le m_0$ 
be given by Lemma~A.4. It is not difficult to see that
$\sup_{z\in \Lambda\setminus  B_1} \|Df^{-n(z)}_{/F(z)}\|'<\mu<1$
(for the real norm). Let
$C=\sup_{z\in \Lambda}\{\|Df^{-j}_{/F(z)}\|'\, \,  \, 1\le j\le m_0\}$. 
Take $k>m_0$ such that
$\mu^k C<\lambda_2< \lambda_1$. Also notice that there exists $r>m_0$ such that if
$z\in \Lambda\setminus  B_1$ and $f^{-j}(z)\in B_1$ for all $n(z)\le j\le n$ with $n\ge r$
then $\|Df^{n}_{/F(z)}\|'<\lambda_2/C$. 

We will show the
lemma for  $m_1= (k+1)r$:  For $z\in \Lambda \setminus B_1$, set
$n_0=0$, $n_1=n(z)$, $n_2=n(f^{n_1}(z)), \ldots$. It follows that for some $i\ge 2$ we
have $m_1=n_1+n_2+...+n_i+s$ with $0\le s\le m_0$. If $s=0$ there is
nothing to do. Consider $j_0=0,j_1, j_2,...,j_\ell$ such that
$f^{n_0+n_1+...+n_{j_s}}(z)\notin B_1$. If $\ell>k$ the result
follows. If $\ell\le k$ then there is $j_s$ such that
$n_{j_{s+1}}-n_{j_s}>r$, and we also conclude. Since we proved  the real claim
for $\lambda_2 < \lambda_1$, there is a small complex
neighbourhood $W_1$  so that the lemma holds for $W_1$
and $\lambda_1$.
\qed
\enddemo

Assume that $\Lambda$ contains only one, say $F$-nonhyperbolic fixed point $q$.
(The general case follows in a very similar way, in particular Lemmas~A.4--5 can
be adapted to the situation where $\NN$ contains more than one point.)
Using $m_1$ from Lemma A.5, we define a second global norm  by:
$$
\|w\|''_{z}=\sum_{j=0}^{m_1-1}\|Df^{-j}w\|'_{f^{-j}z }\, , 
\quad z \in  \Lambda^\complex \,  .\tag{A.8}
$$
More generally, we shall use auxiliary norms
$
\|w\|''_{i, z}=\sum_{j=0}^{i-1}\|Df^{-j}w\|'_{f^{-j}z }
$
for $1\le i \le m_1$.
(In particular $\|w\|''_{1, z} =\|w\|'_{ z}$.) If $z\in W_1$ then by Lemma~ A.5, 
$
\|Df^{-1} _{/F(z)}\|''_{ z}< 1
$.

\medskip

{\bf The global adapted norm and the global adapted semi-norms}

\smallskip

We have two global complex norms
$\|\cdot\|'$ and $\|\cdot\|''$: one which is adapted 
in an admissible neighbourhood $\widetilde B$ of $q$; the other one adapted
outside a neighbourhood $B_1$ of $q$, and we may assume that $B_1$ is
much smaller than $B$.  We must glue these two norms in order to
get a global adapted complex norm $\|\cdot \|_z$  over $\widetilde \Lambda
\subset \Lambda^\complex$.

This can be done as in [Cr, pp. 1112--1114], using the
auxiliary norms $\|w\|'_{-i, z}$,  and $\|w\|^{', \ell}_{-i, z}$ defined there,
because we have the equivalent of (19) and (17)\& (20) in Crovisier's
paper:
$$
\eqalign
{
&\|w\|'_{z} < \|Df w\|'_{f(z)}= \| Df w\|''_{1,  f(z)}<
\| Df^2 w\|''_{2,  f^2(z)}<
\cdots <  \| Df^{m_1} w\|''_{ f^{m_1} z}\, ,\cr
&\qquad\qquad\qquad \qquad\qquad\qquad \forall z \in \widetilde B  \, , 0\ne w \in F(z)\, ;}\tag{A.9}
$$
if $0\ne w \in F(z)$,  and $n-1\ge m_1$, we get from Lemma A.5
and \thetag{A.7}:
$$
\|Df^{n} w\|'_{ z} > 4 \| w \|'_{ z}\, , \forall z
\hbox{ with } f^{n-1} z \notin B_1 \, ;\tag{A.10}
$$
finally, using Lemma~A.5 again,
there is $m_2 \gg m_1$ so that  if $0\ne w \in F(z)$
$$
\|Df^{m_2} w\|'_{ z} >  Cm_1 \|w\|'_z\ge 4 \| w \|''_{ z}\, , \forall z 
\hbox{ with } f^{m_2-1} z \notin B_1 \, . \tag{A.11}
$$
(Note that Lemma~A.5 gives us better control than what is available in [Cr]
so that his construction can in fact be simplified.)

\smallskip
Finally, the two adapted semi-norms are just
$$\|w\|_{E,z}=\| P^E_z w \|_z \, , \quad \|w\|_{F,z}=\| P^F_z w \|_z 
\, , \quad z \in \Lambda^\complex \, .\tag{A.12}
$$ 
The smoothness property \thetag{4.8} follows from the
fact that $z \mapsto \|\cdot\|_z$ is $C^\infty$, combined with the Hoelder smoothness of
the foliations and $C^\infty$ property of leaves.
\qed

\enddocument